\documentclass[12pt]{article}

\usepackage[english]{babel} 		
\usepackage{amsmath}
\usepackage{array}
\usepackage{amssymb}
\usepackage{enumerate,epsfig}		
\usepackage{graphicx}   		
\usepackage{color}			
\usepackage{hyperref}

\newtheorem{Th}{Theorem}[section]
\newtheorem{Prop}[Th]{Proposition}
\newtheorem{Lemma}[Th]{Lemma}

\newtheorem{Rk}[Th]{Remark}

\makeatletter
\@addtoreset{equation}{section}
\makeatother

\newcommand{\cqfd}{{\unskip\kern 6pt\penalty 500
\raise -2pt\hbox{\vrule\vbox to 6pt{\hrule width 6pt
\vfill\hrule}\vrule}\par}}

\def\del{\partial}			
\def\eps{\varepsilon} 			
\def\N{{\mathbb N}} 			
\def\Z{{\mathbb Z}} 			
\def\r{{\mathbb R}} 			

\def\1{{1\hspace{-0.9mm}{\rm l}}}	
\def\ds{\displaystyle}			

\def\tw{\mathop{{\rightharpoonup}}\limits} 

\def\ssi{\quad\Longleftrightarrow\quad}	

\newcommand{\ft}[1]		
{\left(\begin{array}{c} #1 \end{array} \right)}
\newcommand{\coord}[2]		
{\ft{#1\\#2}}
\newcommand{\coords}[3]		
{\ft{#1\\#2\\#3}}
\newcommand{\ligne}[3]	
{\left(\begin{array}{ccccc} #1 & #2 & \dots & #3 \end{array} \right)}
\newcommand{\colonne}[3]	
{\left(\begin{array}{c} #1\\#2\\ \vdots\\#3 \end{array} \right)}
\newcommand{\colonnehaute}[3]	
{\left(\begin{array}{c} #1\\#2\\ \vdots\\\vdots\\#3 \end{array} \right)}
\newcommand{\deter}[4]		
{\left|\begin{array}{cc} #1 & #2\\#3 & #4 \end{array}
\right|}
\newcommand{\mat}[4]		
{\left(\begin{array}{cc} #1 & #2\\#3 & #4 \end{array}
\right)}
\newcommand{\determ}[9]		
{\left|\begin{array}{ccc} #1 & #2 & #3\\#4 & #5
& #6\\ #7 & #8 & #9\end{array}
\right|}
\newcommand{\matr}[9]		
{\left(\begin{array}{ccc} #1 & #2 & #3\\#4 & #5
& #6\\ #7 & #8 & #9 \end{array}
\right)}
\newcommand{\determinant}[9]	
{\left|\begin{array}{cccc} #1 & #2 &\dots &  #3\\#4
& #5
& \dots & #6\\ #7 & #8 & \dots & #9\end{array}
\right|}
\newcommand{\matrice}[9]	
{\left(\begin{array}{cccc} #1 & #2 &\dots &  #3\\#4
& #5
& \dots & #6\\ \vdots & \vdots & \ddots & \vdots\\
#7 & #8 & \dots & #9\end{array}
\right)}
\newcommand{\matricelongue}[9]	
{\left(\begin{array}{ccccc} #1 & #2 &\dots & \dots &  #3\\#4
& #5
& \dots & \dots & #6\\ \vdots & \vdots  & & & \vdots\\
#7 & #8 & \dots & \dots & #9\end{array}
\right)}
\newcommand{\matricehaute}[9]	
{\left(\begin{array}{cccc} #1 & #2 &\dots &  #3\\#4
& #5
& \dots & #6\\ \vdots & \vdots &  & \vdots\\\vdots & \vdots &  &
\vdots\\
#7 & #8 & \dots & #9\end{array}
\right)}

\begin{document}
\title{A model for the evolution of traffic jams in multi-lane} 

\author{F. Berthelin \& D. Broizat \vspace{0.5cm} \\ 
Laboratoire J. A. Dieudonn\'e, UMR 6621 CNRS,\\
        Universit\'e de Nice, Parc Valrose,\\
        06108 Nice cedex 2, France, \\
	Email: Florent.Berthelin@unice.fr, Damien.Broizat@unice.fr}

\date{\today}

\vskip 0.4cm

\maketitle

\abstract{
In \cite{berthelin_degond_delitala_rascle},  Berthelin, Degond, Delitala and Rascle introduced a traffic flow model describing the formation and the dynamics of traffic jams. This model consists of
a Pressureless Gas Dynamics system under a maximal constraint on the density and is derived through a singular limit of the Aw-Rascle model.
In the present paper we propose an improvement of this model by allowing the road to be multi-lane piecewise.
The idea is to use the maximal constraint to model the number of lanes.
We also add in the model a parameter $\alpha$ which model the various speed limitations according to the number of lanes. 
We present the dynamical behaviour of clusters (traffic jams) and by approximation with such solutions,
we obtain an existence result of weak solutions for any initial data.
}



\vskip 0.4cm

\newpage

\tableofcontents

\section{Introduction}
\label{intro}

Classical models of traffic are splitted into three main categories:
particle models (or ``car-following'' models) \cite{gazis_herman_rothery, 
bando_hesebe_nakayama_shibata_sugiyama},
kinetic models \cite{prigogine,prigogine_herman,nelson,klar_wegener},
and fluid dynamical models \cite{lighthill_whitham,payne,payne_bis,aw_rascle,zhang,colombo,helbing}.
Obviously, these models are related; for example in \cite{aw_klar_materne_rascle},
a fluid model is derived from a particle model. See also \cite{klar_kuhne_wegener}
and for recent review on this topic, see \cite{bellomo_dogbe} and \cite{helbing-review}.
Here, we are interested in the third approach, which describes the evolution of 
macroscopic variables (like density, velocity, flow) in space and time.
Let us recall briefly the history of such models.

The simplest fluid models of traffic are based on the single conservation law
$$\del_t n+\del_x f(n)=0,$$
where $n=n(t,x)$ is the density of vehicles and $f(n)$ the associated flow. 
This model only assumes the conservation of the number of cars.
Such models are called ``first order'' models, and the first one is due to Lighthill and Whitham 
\cite{lighthill_whitham} and Richards \cite{richards}.\\
\noindent 
If we take the flux $f(n)=nu$
with  $u=u(t,x)$ the  velocity of the cars,
we add a second equation of equilibrium related to the conservation of momentum.
This approach starts with the Payne-Whitham model \cite{payne,payne_bis}.\\
\noindent 
But the analogy fluid-vehicles is not really convincing: in fact, in the paper \cite{daganzo}, Daganzo shown 
the limits of this analogy, exhibiting absurdities which are implied by classical second-order models,
for example, vehicles going backwards.
To rehabilitate these models, Aw and Rascle proposed in \cite{aw_rascle} a new one which corrects
the deficiencies pointed out by Daganzo. In particular, the density and velocity remain nonnegative.\\
\noindent 
The Aw-Rascle model is given by
$$\left\{\begin{array}{llll}\label{AR_model}
\del_t n+\del_x (nu) = 0,\\
(\del_t +u\del_x)(u+p(n)) =  0,
\end{array}\right.$$
or in the conservative form
$$\left\{\begin{array}{llll}\label{AR_conservative_model}
\del_t n+\del_x (nu)  = 0,\\
\del_t(n(u+p(n))) +\del_x(nu(u+p(n))) =  0,
\end{array}\right.$$
where $p(n)\sim n^{\gamma}$ is the velocity offset, which bears analogies with the pressure in fluid dynamics.\\
\noindent 
In fact, this model can be derived from a microscopic ``car-following'' model, 
as it has been shown in \cite{aw_klar_materne_rascle}.
But even the Aw-Rascle model exhibits some unphysical feature, namely the non-propagation of the upper bound 
of the density $n$, making a constraint such that $n\leq n^*$ impossible (where $n^*$ stands for a maximal 
density of vehicles). \\
\noindent 
Some constraints models have been developed these last years in order to impose such bounds
in hyperbolic models.
See \cite{bouchut_brenier}, \cite{bertheli}, \cite{berthelin_bouchut} for the first results
of this topic and \cite{berthelin_num} for a numerical version of this kind of  problem.\\
\noindent 
That is why recently, Berthelin, Degond, Delitala and Rascle \cite{berthelin_degond_delitala_rascle} proposed a new second-order
model, which aim is to allow to preserve the density constraint $n\leq n^*$ at any time. The main ideas are:
\begin{itemize}
 \item[$\bullet$] modifying the Aw-Rascle model, changing the velocity offset into 
$$p(n)=\left(\frac{1}{n}-\frac{1}{n^*}\right)^{-\gamma},\qquad n<n^*,$$
thus $p(n)$ is increasing and tends to infinity when $n\rightarrow n^*;$
 \item[$\bullet$] rescaling this modified Aw-Rascle model (changing $p(n)$ into $\eps p(n_\eps)$) and taking the
formal limit when $\eps\rightarrow 0^+$.
\end{itemize}
This process leads to a limit system on $(n,u)$ which corresponds to the \textit{Pressureless Gas Dynamics system}:
$$\left\{\begin{array}{llll}\label{PGD_model}
\del_t n+\del_x (nu) = 0,\\
\del_t (nu)+\del_x(nu^2) =  0,
\end{array}\right.$$ 
in areas where $n<n^*$. But a new quantity appears, due to the singularity of the velocity offset in $n=n^*$.
In fact, denoting by $\overline{p}(t,x)$ the formal limit of $\eps p(n_{\eps})(t,x)$ when $\eps\rightarrow 0^+$,
we may have $\overline{p}$ non zero and finite at a point $(t,x)$ such that $n(t,x)=n^*$. Thus, the function 
$\overline{p}$ turns out to be a Lagrangian multiplier of the constraint $n\leq n^*$.
Finally, we obtain 
the \textit{Constrained Pressureless Gas Dynamics} (designed as CPGD) model:
$$\left\{\begin{array}{llll}\label{CPGD_model_conservatif}
\del_t n+\del_x (nu) = 0,\\
\del_t(n(u+\overline{p})) +\del_x(nu(u+\overline{p})) =  0,\\
0 \leq n \leq n^* \, , \quad \overline{p} \geq 0 \, , \quad 
(n^*- n) \overline{p} = 0. 
\end{array}\right.$$
The term $\overline{p}$ represents the speed capability which is not used if the road is blocked and 
that the cars in front imposes a speed smaller than that desired.
We refer to \cite{berthelin_degond_delitala_rascle} for more details on the derivation of the CPGD system in the case 
of the maximal density $n^*$ being constant.
The case where $n^*$ depends on the velocity ($n^*=n^*(u)$) is more realistic (taking into account the fact that the maximal 
number of cars is smaller as the velocity is great, for safety reasons) and is treated
in \cite{berthelin_moutari}.
In \cite{degond_delitala}, a numerical treatment of traffic jam is done.

In this paper, we propose another type of improvement based on the following idea:
the idea is to use the maximal constraint to model the number of lanes.
The constraint $n^*$ will depend on the number of lanes in the portion of the road.
Indeed, in a two-lane portion,
$n^*$ can be twice greater than it is in a one-lane portion of the road.
This idea simplifies the model dramatically and we no longer need to consider as 
many equations of lanes which makes the modeling and use much simpler while reporting the same phenomenon.\\
\noindent 
The new emerging behaviors which are obtained compared to previous models 
\cite{berthelin_degond_delitala_rascle}  and \cite{berthelin_moutari} are the following:
\begin{itemize}
 \item possibility for  cars to accelerate (when the road widens) and then
to change their maximum/wanted velocity,
 \item creation of a void area in a jam (the acceleration of the leading car is not necessarily followed if there is not a sufficient reserve of speed),
 \item this point represents also an approach to model some kind of stop and go waves which is new in such model,
 \item and of course, the multi-lane approach.
\end{itemize} 
The paper is organized as follows: in the next section, we make a modification of the CPGD system to model traffic jams in 
multi-lane. In section 3, we present the dynamics of jams. By approximation with such data, it is
used in section 4 to prove the existence of weak solutions for any initial data.

\section{The ML-CPGD model}
We consider a piecewise constant maximal density of vehicles, given by
\[
n^*(x) = \sum_{j=0}^{M} n^*_j \1_{]r_j,r_{j+1}[}(x)
\]
where
\[
n_j^*\in \{1, 2\},\quad (r_j)_{1\leq j\leq M} \textrm{
an increasing
sequence of real numbers}, 
\]
\[
r_0=-\infty,\quad r_{M+1}=+\infty.
\]
It means that we set on a road with one or two lanes, the road transitions (change of number of lanes)
being at points $(r_j)_{1\leq j\leq M}$. On a one-lane section, the maximal
density is one (in
view
of simplification), whereas on a two-lane section, the maximal allowed density
is two. 
It is the first improvement of our model: the constraint density changes with $x$ to model the fact that there is one or two lanes.
Evolution equations are given by the \textit{Multi-lane Constrained Pressureless Gas Dynamics} 
system (designed by ML-CPGD), whose conservative form is
\begin{eqnarray}
& &       \partial_t n +  \partial_x (n u) = 0 \, ,  \label{CPGD1}\\
& &      \partial_t (n(u+p)I_{\alpha})+  \partial_x (n u(u + p )I_{\alpha}) = 0
\, , \label{CPGD2}\\
& & 0 \leq n \leq n^*(x) \, , \quad \, u\geq 0 \,, \quad p \geq 0 \, , \quad 
(n^*(x) - n) p = 0 \, 
  \label{CPGD_cons1},
\end{eqnarray}
where the function $I_{\alpha}=I_{\alpha}(x)$ is defined by
$$ I_{\alpha}(x)=\left\{\begin{array}{lll}
1 & \mbox{if} & n^{*}(x)=1,\\
1/\alpha & \mbox{if} & n^{*}(x)=2.
  \end{array}\right.$$
The number $\alpha\geq 1$ stands for the rate between two-lane velocities
and one-lane velocities.
Thus a single car (we mean a car not into a jam) with speed $u$ on a one-lane road will pass to the speed $\alpha u$ on a two-lane road. 
This represents the fact that on a two-lane section, the average velocity is
higher than on a one-lane (on a highway, you drive faster than on a road
even if you are alone). The preferred velocity depends on the road width
according to $\alpha$.
It can also be understood as the speed limitation on the various kind of roads. This is the second improvement of our model.
It only act on the second equation since it is the momentum quantity which has to be changed and not the conservation of the number of cars 
(first equation). 
\vspace{5mm}

Of course, this model can be extended to case with  three-lane, four-lane portion...
In the case of three lanes, $n_j^*\in \{1, 2, 3\}$ and $I_\alpha$ is replaced by 
$\displaystyle I_{\alpha, \beta}(x)=\left\{\begin{array}{lll}
1 & \mbox{if} & n^{*}(x)=1,\\
1/\alpha & \mbox{if} & n^{*}(x)=2,\\
1/\beta & \mbox{if} & n^{*}(x)=3,
  \end{array}\right.$
with $\beta \geq \alpha \geq 1$,
$\alpha$ being the rapport of speed between one and two lanes and
$\beta/\alpha$ the rapport between three and two lanes.

\section{Clusters dynamics}\label{bouchons_collants}

In this section, we present some particular solutions $(n,u,p)$ of (\ref{CPGD1})-(\ref{CPGD_cons1})
which are clusters solutions. For these functions, $n=n(t,x)$ take as only values $0$ and $n^{*}(x)$.\\
\noindent 
In some sense, they are an extension of sticky particles of
\cite{brenier_grenier,e_rykov_sinai} playing a crucial role in the proof of existence of solutions
for constraint models.
They have been introduced in \cite{bouchut_brenier} and used with various dynamics in
\cite{bertheli,berthelin_bouchut,berthelin_degond_delitala_rascle,berthelin_moutari}.\\
\noindent 
Let us consider the density $n(t,x)$, the flux $n(t,x)u(t,x)$  and the pressure
$n(t,x) p(t,x)$ given respectively by
\begin{align}
&n(t,x)= n^{*}(x)\sum_{i=1}^{N}\1_{a_{i}(t)<x<b_{i}(t)},\label{defn} \\
&n(t,x) u(t,x)=n^{*}(x)\sum_{i=1}^{N}u_{i}\1_{a_{i}(t)<x<b_{i}(t)},\label{defnu}\\
&n(t,x) p(t,x)=n^{*}(x)\sum_{i=1}^{N}p_{i}\1_{a_{i}(t)<x<b_{i}(t)},\label{defnp}
\end{align} 
with
\[N\in\N^{*},\qquad u_{i}\geq 0,\qquad p_{i}\geq 0,\]
as long as there is no collision and no change of $n^*(x)$. 
That is to say
\[a_1(t)<b_1(t)\leq a_2(t)<b_2(t)<\ldots\leq a_N(t)<b_{N}(t)\]
and the number of blocks $N$ is constant until there is a shock or a change of width (thus we have $N=N(t)$).\\
\noindent 
This type of piecewise constant solution writes as a superposition of blocks with $(n,u,p)=(n^{*},u_{i},p_{i})$ constant.
Each block evolves according to the interactions with the other blocks and the changes of width.\\
\noindent 
We have to explain three dynamics:
\begin{itemize}
 \item What happens when two blocks collide ? (how to describe a shock) 
 \item What happens when the road narrows ($n^*(x)$ was 2 and becomes 1) ?
 \item What happens when the road widens ($n^*(x)$ was 1 and becomes 2) ?
\end{itemize}

\noindent First, let us present some technical properties that will be used in the various cases.

\begin{Lemma}
Let be $s,\sigma\in [0,+\infty[$, 
$a,b\in C^1(]inf(s,\sigma),sup(s,\sigma)[)$
and $\varphi\in\mathcal{D}(]0,+\infty[\times\r)$. We set
$$J(s,\sigma,a,b,u):=\int_{s}^{\sigma}\int_{a(t)}^{b(t)}(\del_t\varphi(t,
x)+u(t)\del_x\varphi(t,x))dxdt.$$ 
Then we get \begin{eqnarray}\label{calcul}
J(s,\sigma,a,b,u)&=&\int_{a(\sigma)}^{b(\sigma)} \varphi(\sigma,x) \,dx
-  \int_{a(s)}^{b(s)} \varphi(s,x) \,dx \\
& &+\int_{s}^\sigma\varphi(t,b(t)) \, (u(t)-b'(t)) \,dt
          +\int_{s}^\sigma\varphi(t,a(t)) \, (a'(t)-u(t)) \,dt.  \nonumber
\end{eqnarray}
\end{Lemma}

\noindent\underline{Proof}: 
We have
\[\ds 
\frac{d}{dt} \left[ 
\int_{a(t)}^{b(t)} \varphi(t,x) \,dx
\right]
=\int_{a(t)}^{b(t)} \partial_t \varphi(t,x) \,dx +\varphi(t,b(t)) \, b'(t)
           -\varphi(t,a(t)) \, a'(t),\]
thus
\[ \begin{split}
\int_{s}^\sigma \int_{a(t)}^{b(t)} \partial_t \varphi(t,x) \,dx \,dt 
 & = 
\int_{a(\sigma)}^{b(\sigma)} \varphi(\sigma,x) \,dx
-  \int_{a(s)}^{b(s)} \varphi(s,x) \,dx \\
 & \quad -\int_{s}^\sigma\varphi(t,b(t)) \, b'(t) \,dt
          +\int_{s}^\sigma\varphi(t,a(t)) \, a'(t)\,dt.\end{split}
\] 
Moreover
$\ds 
\int_{s}^\sigma \int_{a(t)}^{b(t)} \partial_x \varphi(t,x) \,dx \,dt 
=
\int_{s}^{\sigma} \varphi(t,b(t)) \,dt
-  \int_{s}^{\sigma} \varphi(t,a(t)) \,dt$ and
the result follows.
\cqfd

\begin{Rk} We notice that
\begin{equation}\label{symetrie 1 de J}
J(\sigma,s,a,b,u)=-J(s,\sigma,a,b,u), 
\end{equation}
\begin{equation}\label{symetrie 2 de J}
J(s,\sigma,b,a,u)=-J(s,\sigma,a,b,u).
\end{equation}
If we have $a'=b'=u$, then 
\begin{equation}\label{simplification de J}
J(s,\sigma,a,b,u)=\int_{a(\sigma)}^{b(\sigma)} \varphi(\sigma,x) \,dx
-  \int_{a(s)}^{b(s)} \varphi(s,x) \,dx.
\end{equation}
If we have $a'=u$ and $c$ is constant, then 
\begin{equation}\label{simplification 2 de J}
J(s,\sigma,a,c,u)=\int_{a(\sigma)}^{c} \varphi(\sigma,x) \,dx
-  \int_{a(s)}^{c} \varphi(s,x) \,dx+\int_{s}^{\sigma}\varphi(t,c)u(t)dt.
\end{equation}
\end{Rk}

\begin{Lemma}
 We have the following formulas:\\
\noindent If $a'=b'=c'=u$, then
\begin{equation}\label{formule 1 sur J}
\begin{split}
J(s,\sigma,a,b,u)+J(\sigma,\tau,a,c,u) & 
  =\ds-\int_{a(s)}^{b(s)} \varphi(s,x)
\,dx \\ & \quad +\int_{c(\sigma)}^{b(\sigma)}
\varphi(\sigma,x) \,dx + \int_{a(\tau)}^{c(\tau)}\varphi(\tau,x) \,dx.	
\end{split}
\end{equation}
If $a'=b'=u$ and $c=b(\sigma)=a(\tau)$, then
\begin{equation}\label{formule 2 sur J}
J(s,\sigma,a,b,u)+J(\sigma,\tau,a,c,u)=-\int_{a(s)}^{b(s)}\varphi(s,x)dx+\int_{
\sigma}^{\tau}u(t)\varphi(t,c)dt.
\end{equation}
\end{Lemma}

\noindent\underline{Proof}:
We have
\begin{eqnarray}
& & \hspace{-10mm} J(s,\sigma,a,b,u)+J(\sigma,\tau,a,c,u)\\  
& =  &  \ds\int_{a(\sigma)}^{b(\sigma)}
\varphi(\sigma,x) \,dx-\int_{a(s)}^{b(s)} \varphi(s,x) \,dx\\
& & +\ds\int_{s}^\sigma\varphi(t,b(t)) \, (u(t)-b'(t)) \,dt
+\int_{s}^\sigma\varphi(t,a(t)) \, (a'(t)-u(t)) \,dt\\
& &+  \ds\int_{a(\tau)}^{c(\tau)}
\varphi(\tau,x) \,dx-\int_{a(\sigma)}^{c(\sigma)} \varphi(\sigma,x) \,dx\\
& &+  \ds\int_{\sigma}^{\tau}\varphi(t,c(t)) \, (u(t)-c'(t)) \,dt
+\int_{\sigma}^{\tau}\varphi(t,a(t)) \, (a'(t)-u(t)) \,dt\\
& = & \ds-\int_{a(s)}^{b(s)} \varphi(s,x) \,dx +
\int_{\sigma}^{\tau}\varphi(t,c(t)) \, (u(t)-c'(t))
\,dt\\
& &+  \ds\int_{c(\sigma)}^{b(\sigma)}
\varphi(\sigma,x) \,dx + 
\int_{a(\tau)}^{c(\tau)}\varphi(\tau,x) \,dx\\
& &+  \ds\int_{s}^\sigma\varphi(t,b(t)) \, (u(t)-b'(t)) \,dt
+\int_{s}^{\tau}\varphi(t,a(t)) \, (a'(t)-u(t)) \,dt.
\end{eqnarray}
Since $a'=b'=u,$ the two last terms vanish and we have
$$\begin{array}{llll}
 J(s,\sigma,a,b,u)+J(\sigma,\tau,a,c,u)
& = & \ds-\int_{a(s)}^{b(s)} \varphi(s,x)
\,dx\\ 
 &  &+
\ds\int_{\sigma}^{\tau}\varphi(t,c(t)) \, (u(t)-c'(t))
\,dt\\
& & + \ds\int_{c(\sigma)}^{b(\sigma)}
\varphi(\sigma,x) \,dx + 
\int_{a(\tau)}^{c(\tau)}\varphi(\tau,x) \,dx.
\end{array}$$
The formulas \eqref{formule 1 sur J} and \eqref{formule 2 sur J} follow.
\cqfd

\subsection{About uniqueness of the dynamics}\label{uniqueness dynamics}

\noindent In order to work with the most realistic solution, it is necessary to impose a certain number of criteria on the dynamics in question.
This discussion also improve the paper \cite{berthelin_degond_delitala_rascle}.

A single block for which $u+p$ stays constant is a solution, for example the function corresponding to the following figure:
\begin{center}
\includegraphics[width=0.6\textwidth]{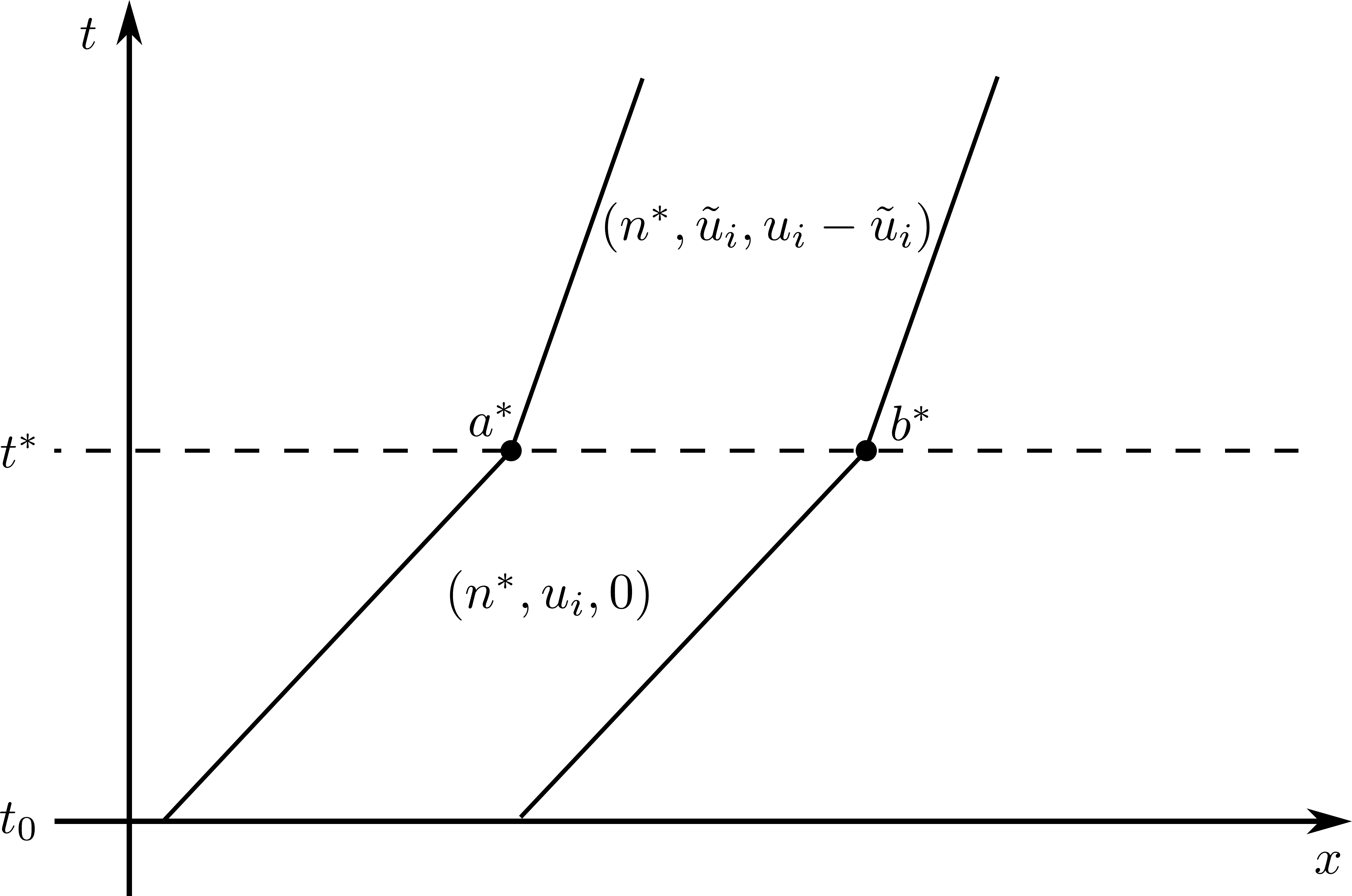}
\end{center}
\begin{Rk}
To understand the meaning of the dynamics,
for every figure, the term $(n,u,p)$ on a zone corresponds to the constant values of the functions on a block. 
\end{Rk}
In fact, in an open subset $\Omega\subset ]0,+\infty[_{t}\times \r_{x}$, where $n^{*}$ is constant, it is very easy to see that the dynamic 
displayed on figure satisfies \eqref{CPGD1}-\eqref{CPGD_cons1}, for any value of $0 \leq \tilde{u}_{i} \leq u_i$.\\
\noindent 
Now, remember that the term $p$ represents the speed capability which is not used if the road is blocked and 
that the cars in front imposes a speed smaller than that desired.
The term is 0 if the density is not $n^*(x)$ since in this case the car can go to its preferred velocity.
Thus there is no reason for a single car to have a nonzero pressure term if there is no one before him.
And the relation $(n^{*}(x)-n)p=0$ do not impose $p=0$ for the first car of the jam.
This is why we  assume that the blocks satisfy the additional constraint:
\begin{equation} \label{cont+}
(n^{*}(x)-n(x^{+}))p=0, 
\end{equation}
in zones where $n^{*}$ is constant.
In this property, we denote by $n(x^{+})$ the limit, if it exists, of $n(y)$ when $y \to x$ with $y>x$.
With this condition, the dynamics of the previous figure is a solution only if
$\tilde u_i=u_i$.

The interpretation is the following :
if the first car of the jam has the opportunity to use its preferred velocity, it uses it and $p$ becomes zero.
If not, $p$ is not necessarily zero. \\
\noindent 
This is why for various blocks sticking one after the other, the constraint (\ref{cont+}) on $p$ gives
the two situations of the above figures.

\begin{center}
\includegraphics[width=0.6\textwidth]{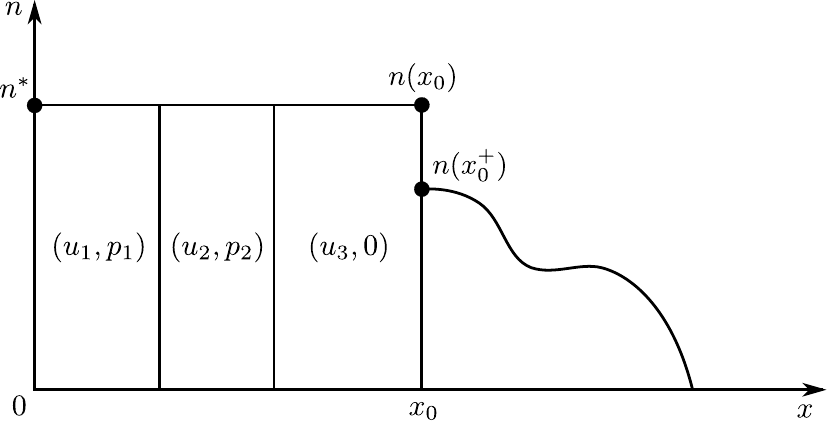}\hspace{10mm}
\includegraphics[width=0.6\textwidth]{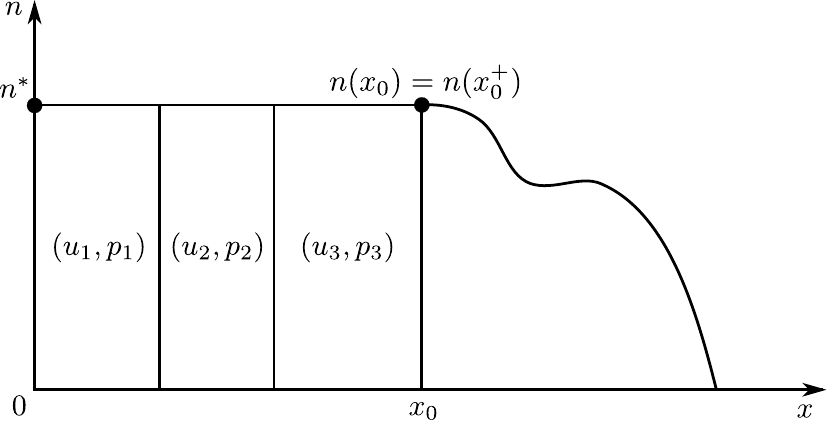}
\end{center}

In fact, an other criteria than (\ref{cont+}) related to the minimization of $p\geq 0$ can be used.
It is clear than choosing $p=0$ minimize the $p$ term in the previous described situations.
When the constraint (\ref{cont+}) cannot be imposed, uniqueness criteria which is natural is the minimization of $p\geq 0$.

\noindent 
We now detail  the various cases that can appear in the dynamic of clusters.

\subsection{Collision between two blocks without change of width} 

In a zone where $n^*(x)=n^*$ is constant, we consider two blocks $(n^{*},u_{l},0)$ and $(n^{*},u_{r},0)$, with $u_{l}>u_{r}$.
Thus, at a time $t^{*}>0$, the left block reaches the right one, and collide with it.
The dynamic is displayed in the following figure.
\begin{center}
\includegraphics[width=0.6\textwidth]{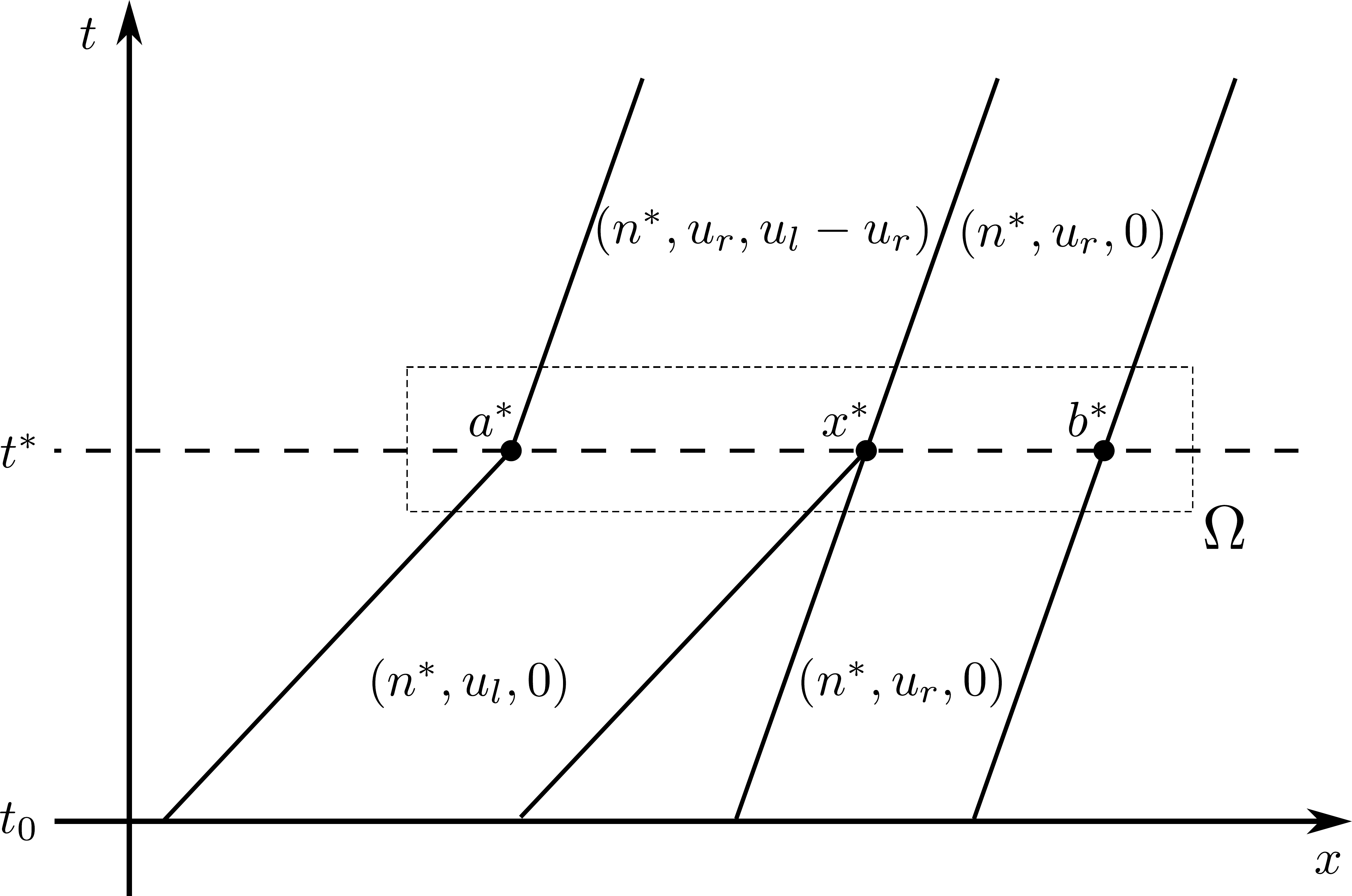}	
\end{center}
The density $n(t,x)$, the flux $n(t,x) u(t,x)$ and the function $p(t,x)$
are locally given respectively by  
\[
n(t,x)= 
\left\{\begin{array}[]{lll}
n^{*}\1_{a_{l}(t)<x<b_{l}(t)}+n^* \1_{a_{r}(t)<x<b_{r}(t)} & \mbox{if} & t<t^*,\\
n^{*} \1_{\tilde a_{l}(t)<x<\tilde b_{l}(t)}+n^*\1_{a_{r}(t)<x<b_{r}(t)} & \mbox{if} & t>t^*,
\end{array}
\right.
\]
\[
n(t,x)u(t,x)= \left\{\begin{array}[]{lll}
n^{*} u_{l} \1_{a_{l}(t)<x<b_{l}(t)}+n^* u_{r}
\1_{a_{r}(t)<x<b_{r}(t)} & \mbox{if} & t<t^*,\\
n^{*}u_{r} \1_{\tilde a_{l}(t)<x<\tilde b_{l}(t)}+n^*u_{r}
\1_{a_{r}(t)<x<b_{r}(t)} & \mbox{if} & t>t^*,
\end{array}\right.
\]
and 
\[
n(t,x)p(t,x)= 
\left\{\begin{array}[]{lll}
0 & \mbox{if} & t<t^*,\\
n^{*}(u_{l}-u_{r}) \1_{\tilde a_{l}(t)<x<\tilde b_{l}(t)} & \mbox{if} & t>t^*,
\end{array}\right.
\]
with the linear functions $a_{l}$, $b_{l}$, $a_{r}$, $b_{r}$, $\tilde a_{l}$, $\tilde b_{l}$ are given by
\[
\frac{d}{dt}a_l(t)=\frac{d}{dt}b_l(t)=u_{l},\quad  a_l(t^*)=a^*,\quad
b_l(t^*)=x^*,
\]
\[
\frac{d}{dt}a_r(t)=\frac{d}{dt}b_r(t)=u_{r},\quad  a_r(t^*)=x^*,\quad
b_r(t^*)=b^*,
\]
\[
\frac{d}{dt}\tilde a_l(t)=\frac{d}{dt}\tilde b_l(t)=u_{r},\quad 
\tilde a_l(t^*)=a^*,\quad
\tilde b_l(t^*)=x^*,
\]
and 
\[u_l>u_r.\] 
The left block obtains the velocity of
the one being immediately on its right when they collide.
We extend this
when more than two blocks collide at a time $t^*$, by  forming
a new block with the velocity of the block on the right of the group.

\begin{Lemma}
The previous dynamic satisfies \eqref{CPGD1}-\eqref{CPGD_cons1}.	
\end{Lemma} 

\noindent\underline{Proof}:
Let $\Omega$ be an open neighborhood of the shock zone (displayed in the previous figure). Then,
we have, for any continuous function $S$ and any test function
$\varphi\in\mathcal{D}(\Omega)$,

$$\begin{array}[]{llll}
& &
\hspace{-15mm}\langle\,\,\del_t(nS(u,p,I_{\alpha}))+\del_x(nuS(u,p,I_{\alpha})),
\varphi\,\,\rangle\\\\
& =
&-\ds\int_0^{+\infty}\!\!\!\int_{\r}n(t,x)S(u(t,x),p(t,x),I_{\alpha}
(x))(\del_t\varphi+u\del_x\varphi)dxdt\\\\
& = & -n^{*}S(u_l,0,I_{\alpha})J(0,t^{*},a_l,b_l,u_l)
\\\\
&  & -n^{*}S(u_{r},u_{l}-u_{r},I_{\alpha})J(t^{*},\infty,\tilde
a_l, \tilde b_{l},u_{r})\\\\
& & - n^{*}S(u_{r},0,I_{\alpha})J(0,\infty, a_{r},b_{r},u_{r})\\\\
& = & \ds\left(-n^{*}S(u_{l},0,I_{\alpha})+n^{*}S(u_{r},u_{l}-u_{r},I_{\alpha})\right)\int_{a^{*}}^{x^{*}}\varphi(t^{*},x)dx.
\end{array}$$
For $S(u,p,I_{\alpha})=1$, we get
$$\langle\,\,\partial_t n +\partial_x (n u),\varphi\,\,\rangle=0.$$
For $S(u,p,I_{\alpha})=(u+p)I_{\alpha}$, we get
$$\begin{array}[]{lll}
\langle\,\,\partial_t (n(u+p)I_{\alpha}) +\partial_x (n
u(u+p)I_{\alpha}),\varphi\,\,\rangle & = & 0.
  \end{array}$$
\cqfd

\subsection{Narrowing of the road without collision}\label{retrecissement}

Let us move to the situation where the road narrows ($n^*(x)$ was
2 and becomes 1). Here, we describe the evolution of a block which undergoes this narrowing. 
The speed will be divided by $\alpha$.\\
\noindent 
The dynamic of the block is exhibited in the following figure.

\begin{center}
\includegraphics[width=0.6\textwidth]{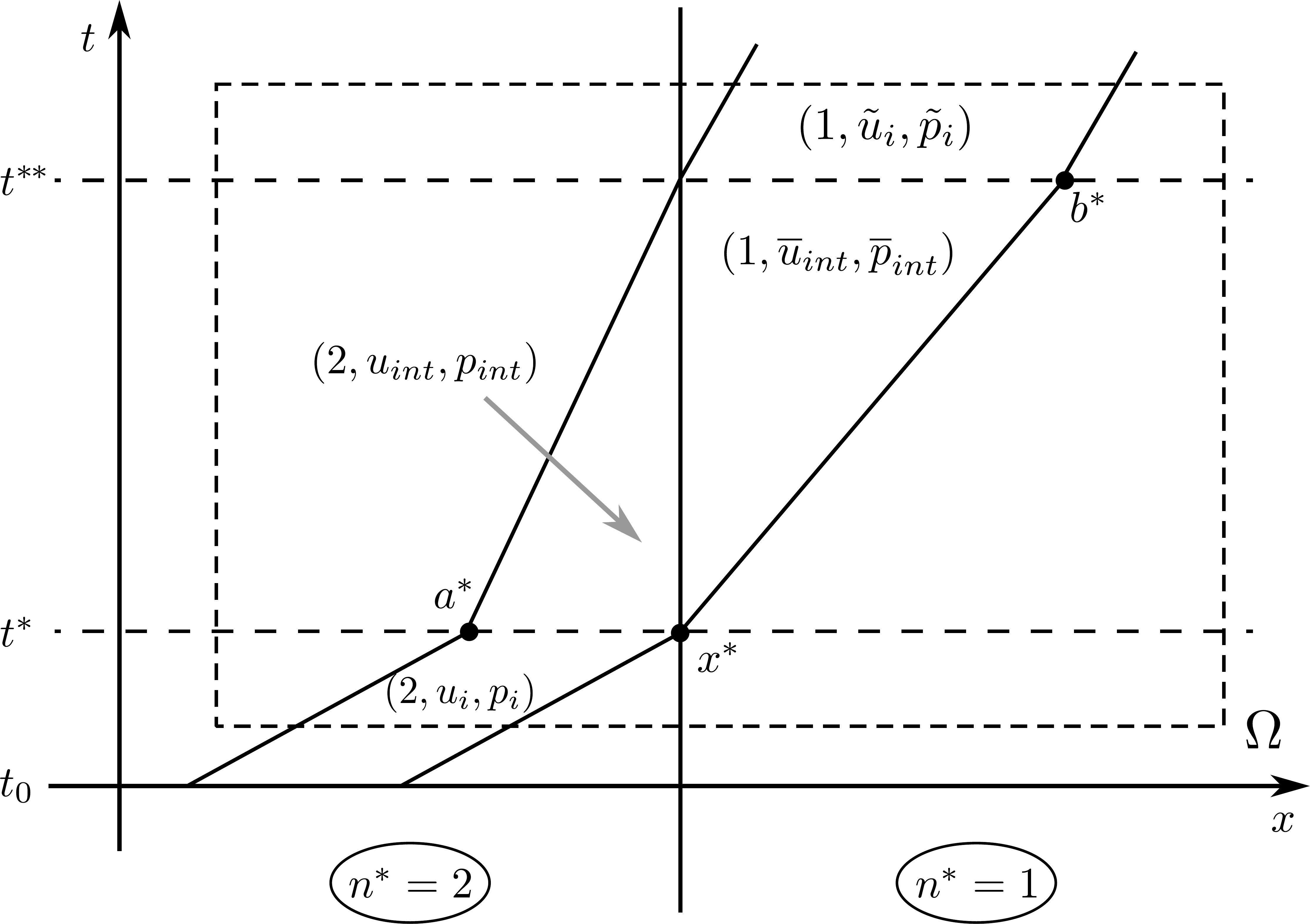}
\end{center}
The density $n(t,x)$, the flux $n(t,x) u(t,x)$ and the functional $p(t,x)$ are
locally given respectively by  
\[
n(t,x)=
\begin{cases}
2\1_{a_{i}(t)<x<b_{i}(t)} & \text{if} \ \ t<t^*,\\
2\1_{a_{int}(t)<x<x^*}+ \1_{x^*<x<\overline b_{int}(t)} & \text{if}
\ \
t^*<t<t^{**},\\
\1_{\tilde a_{i}(t)<x<\tilde b_{i}(t)} &  \text{if} \ \ t>t^{**},
\end{cases}
\]
\[
n(t,x)u(t,x)=
\begin{cases}
 2u_i\1_{a_{i}(t)<x<b_{i}(t)} &  \text{if} \ \ t<t^*,\\
 2u_{int} \1_{a_{int}(t)<x<x^*}+ \overline
u_{int}\1_{x^*<x<\overline b_{int}(t)} & 
\text{if} \
\ t^*<t<t^{**},\\
\tilde u_i \1_{\tilde a_{i}(t)<x<\tilde b_{i}(t)} & \text{if} \ \
t>t^{**},
\end{cases}
\]
and 
\[
n(t,x) p(t,x)=
\begin{cases}
2 p_i \1_{a_{i}(t)<x<b_{i}(t)} & \text{if} \ \ t<t^*,\\
2 p_{int} \1_{a_{int}(t)<x<x^*}+ \overline p_{int}
\1_{x^*<x<\overline
b_{int}(t)} & 
\text{if} \ \ t^*<t<t^{**},\\
\tilde p_i \1_{\tilde a_{i}(t)<x<\tilde b_{i}(t)} &  \text{if} \
\ t>t^{**},
\end{cases}
\]
with

\[
\frac{d}{dt}a_i(t)=\frac{d}{dt}b_i(t)=u_i,\quad  a_i(t^*)=a^*,\quad
b_i(t^*)=x^*,
\]
\[
\frac{d}{dt}a_{int}(t)=u_{int},\quad  a_{int}(t^*)=a^*,\quad
a_{int}(t^{**})=x^*,
\]
\[
\frac{d}{dt}\overline{b}_{int}(t)=\overline{u}_{int},\quad 
\overline{b}_{int}(t^*)=x^*,\quad
\overline{b}_{int}(t^{**})=b^*,
\]
\[
\frac{d}{dt}\tilde a_i(t)=\frac{d}{dt}\tilde b_i(t)=\tilde u_i,\quad 
\tilde a_i(t^{**})=x^*,\quad
\tilde b_i(t^{**})=b^*.
\] 

\begin{Lemma}\label{CNS retrecissement}
The previous dynamic satisfies \eqref{CPGD1}-\eqref{CPGD_cons1} if and only if
\[p_{int}=u_i+p_i-u_{int},\quad
(\overline{u}_{int},\overline{p}_{int})=\left(2 u_{int},\frac{u_{i}+p_{i}}{\alpha}-2 u_{int}\right),\quad 
\tilde p_i=\frac{u_{i}+p_{i}}{\alpha}-\tilde{u}_i,\]
with 
\[0\leq u_{int}\leq \frac{u_i+p_i}{2\alpha},\quad 0\leq \tilde u_i\leq\frac{u_{i}+p_{i}}{\alpha}.\]
\end{Lemma} 

\noindent\underline{Proof}:
We have, for any continuous function $S$ and any function
$\varphi\in\mathcal{D}(\Omega)$,
$$\begin{array}{llll}
& & \hspace{-15mm}\langle\,\,\del_t(nS(u,p,I_{\alpha})+\del_x(nuS(u,p,I_{\alpha})),\varphi\,\,\rangle\\
& =&-\ds\int_0^{+\infty}\!\!\!\int_{\r}n(t,x)S(u(t,x),p(t,x),I_{\alpha}(x))(\del_t\varphi+u\del_x\varphi)dxdt\\\\
& = &-2S(u_i,p_i,1/\alpha)J(0,t^{*},a_i,b_i,u_i)\\\\
& & -2S(u_{int},p_{int},1/\alpha)J(t^*,t^{**},a_{int},x^*,u_{int})\\\\
&  & - S(\overline u_{int},\overline p_{int},1)J(t^{*},t^{**},x^*,\overline b_{int},\overline u_{int})
- S(\tilde u_{i},\tilde p_{i},1)J(t^{**},\infty,\tilde a_i,\tilde b_i,\tilde u_i)\\\\
& =& \ds-2S(u_i,p_i,1/\alpha)\int_{a^*}^{x^*}\varphi(t^*,x)dx\\
& & \ds-2S(u_{int},p_{int},1/\alpha)\left(-\int_{a^*}^{x^*}\varphi(t^*,x)dx+u_{int}\int_{t^*}^{t^{**}}\varphi(t,x^*)dt\right)\\
& & \ds- S(\overline u_{int},\overline p_{int},1)\left(\int_{x^*}^{b^*}\varphi(t^{**},x)dx
-\overline u_{int}\int_{t^*}^{t^{**}}\varphi(t,x^*)dt\right)\\
& & +\ds S(\tilde u_{i},\tilde p_{i},1)\int_{x^*}^{b^*}\varphi(t^{**},x)dx,
\end{array}$$
thus we have
$$\begin{array}{llll}
& & \hspace{-15mm}\del_t(nS(u,p,I_{\alpha})+\del_x(nuS(u,p,I_{\alpha}))\\\\
& = &  2\left(S(u_{int},p_{int},1/\alpha)-S(u_i,p_i,1/\alpha)\right)\1_{[a^*,x^*]}
(x)\delta(t-t^*)\\\\
&  &+ \left(
S(\tilde u_{i},\tilde p_{i},1)- S(\overline u_{int},\overline
p_{int},1)\right)\1_{[x^*,b^*]}
(x)\delta(t-t^{**})\\\\
&  &+  \left(\overline u_{int} S(\overline u_{int},\overline
p_{int},1)-2 u_{int}S(u_{int},p_{int},1/\alpha)\right)\1_{[t^*,t^{**}]}
(t)\delta(x-x^*).
\end{array}$$
For $S(u,p,I_{\alpha})=1$, we get
$$\partial_t n +\partial_x (n u)=\left(\overline
u_{int}-2u_{int}\right)\1_{[t^*,t^{**}]}
(t)\delta(x-x^*).$$
For $S(u,p,I_{\alpha})=(u+p)I_{\alpha}$, we get
$$\begin{array}{lll}
& &\hspace{-15mm}\partial_t (n(u+p)I_{\alpha}) +\partial_x (n u(u+p)I_{\alpha})\\\\
& = &  
\frac{2}{\alpha}\left(u_{int}+p_{int}-u_i-p_i\right)\1_{[a^*,x^*]}
(x)\delta(t-t^*)\\\\
& &+ \left(
\tilde u_{i}+\tilde p_{i}-\overline u_{int}-\overline
p_{int}\right)\1_{[x^*,b^*]}
(x)\delta(t-t^{**})\\\\
&  & +\left(\overline u_{int} (\overline u_{int}+\overline
p_{int})-2 u_{int}\left(\frac{u_{int}+p_{int}}{\alpha}\right)\right)\1_{[t^*,t^{**}]}
(t)\delta(x-x^*).
\end{array}$$
Therefore, $(n,u,p)$ is a solution of \eqref{CPGD1}-\eqref{CPGD2} if and only if
\[
\left\{
\begin{array}[]{lll}
\overline u_{int} & = & 2 u_{int}\\
u_{int}+p_{int} & = & u_i+p_i\\
\tilde u_{i}+\tilde p_{i} & = & \overline u_{int}+\overline
p_{int}\\
\overline u_{int}+\overline p_{int} & = & \frac{1}{\alpha}(u_{int}+p_{int})
\end{array}\right.\quad\Longleftrightarrow\quad
\left\{
\begin{array}[]{lll}
\overline u_{int} & = & 2 u_{int}\\
u_{int}+p_{int} & = & u_i+p_i\\
\overline u_{int}+\overline p_{int} & = & \frac{1}{\alpha}(u_{i}+p_{i})\\
\tilde u_{i}+\tilde p_{i} & = & \frac{1}{\alpha}(u_{i}+p_{i})
\end{array}\right..
\]
Since $u_{int},p_{int}$,$\overline{u}_{int}$ and $\overline{p}_{int}$ are nonnegative, it concludes the proof of lemma.
\cqfd

Now, we can find the dynamics governing a single block $(n^{*},u_{i},0)$ which undergoes a narrowing of the road:
according to subsection \ref{uniqueness dynamics} and additional constraint (\ref{cont+}), we have 
\[p_{i}=0,\qquad \tilde{p}_{i}=0, \qquad \overline{p}_{int}=0,\]
in the relations of lemma \ref{CNS retrecissement}, which leads to 
\[p_{int}=u_i-u_{int},\quad
\overline{u}_{int}=2 u_{int},\quad  u_{int}= \frac{u_i}{2\alpha}, \quad
\tilde{u}_{i}=\frac{u_{i}}{\alpha}.\]
Finally, the only dynamics compatible with (\ref{cont+})  for a narrowing is:
\begin{center}
\includegraphics[width=0.6\textwidth]{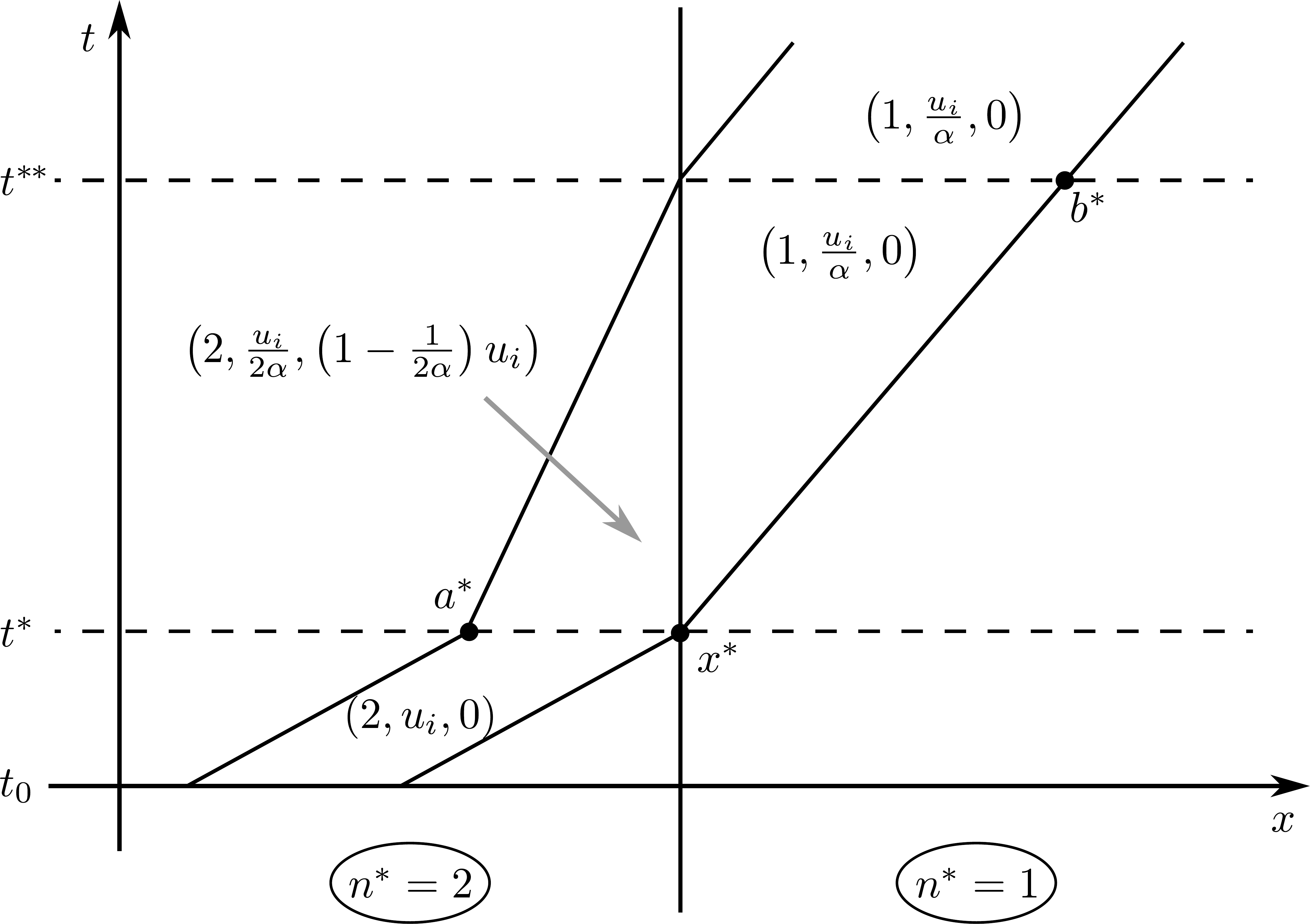}
\end{center}
\begin{Rk}
In this situation, there is a backward propagation of the queue: before the narrowing of the road, the cluster size was $x^{*}-a^{*}$ and during the intermediate state, it increases linearly up to $b^{*}-x^{*}$, with finite speed $\ds\frac{u_{i}}{2\alpha}=\frac{d}{dt}(\overline{b}_{int}(t)-a_{int}(t))$. Concerning the velocity $u$ into the cluster, the information travels at infinite speed and the velocity at the end of the block pass instantly from $u$ to $u/(2\alpha)$.
\end{Rk}
\subsection{Enlargement of the road without collision}\label{elargissement}

Now we explain what happens for a block when the road widens ($n^*(x)$ was 1 
and becomes 2). 

In fact, the block (which comes with $n=n^{*}=1$) becomes a block with $n=n^{*}=2$,
but its speed will be multiplied by the parameter $\alpha$.
The dynamic is exhibited hereafter:

\begin{center}
\includegraphics[width=0.6\textwidth]{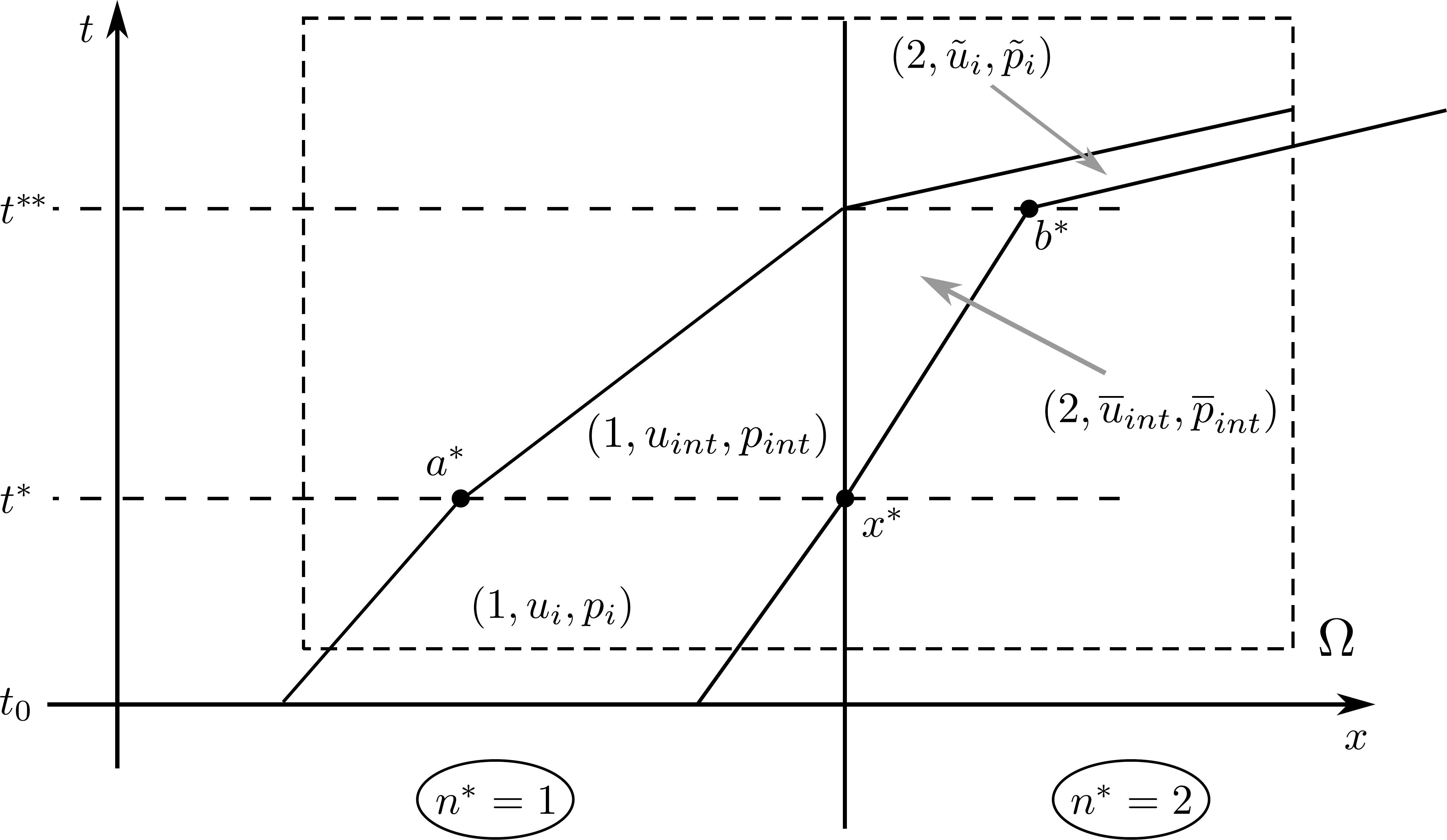}
\end{center}
The density $n(t,x)$, the flux $n(t,x) u(t,x)$ and the functional $p(t,x)$ are
locally given respectively by  
\[
n(t,x)=
\begin{cases}
\1_{a_{i}(t)<x<b_{i}(t)} & \text{if} \ \ t<t^*,\\
\1_{a_{int}(t)<x<x^*}+2\1_{x^*<x<\overline b_{int}(t)} & \text{if} \ \
t^*<t<t^{**},\\
2 \1_{\tilde a_{i}(t)<x<\tilde b_{i}(t)} &  \text{if} \ \ t>t^{**},
\end{cases}
\]
\[
n(t,x)u(t,x)=
\begin{cases}
 u_i\1_{a_{i}(t)<x<b_{i}(t)} &  \text{if} \ \ t<t^*,\\
 u_{int} \1_{a_{int}(t)<x<x^*}+ 2 \overline
u_{int}\1_{x^*<x<\overline b_{int}(t)} & 
\text{if} \
\ t^*<t<t^{**},\\
2\tilde u_i \1_{\tilde a_{i}(t)<x<\tilde b_{i}(t)} & \text{if} \ \
t>t^{**},
\end{cases}
\]
and 
\[
n(t,x) p(t,x)=
\begin{cases}
 p_i \1_{a_{i}(t)<x<b_{i}(t)} & \text{if} \ \ t<t^*,\\
 p_{int} \1_{a_{int}(t)<x<x^*}+ 2\overline p_{int} \1_{x^*<x<\overline
b_{int}(t)} & 
\text{if} \ \ t^*<t<t^{**},\\
2\tilde p_i \1_{\tilde a_{i}(t)<x<\tilde b_{i}(t)} &  \text{if} \
\ t>t^{**},
\end{cases}
\]
with

\[
\frac{d}{dt}a_i(t)=\frac{d}{dt}b_i(t)=u_i,\quad  a_i(t^*)=a^*,\quad
b_i(t^*)=x^*,
\]
\[
\frac{d}{dt}a_{int}(t)=u_{int},\quad  a_{int}(t^*)=a^*,\quad
a_{int}(t^{**})=x^*,
\]
\[
\frac{d}{dt}\overline{b}_{int}(t)=\overline{u}_{int},\quad 
\overline{b}_{int}(t^*)=x^*,\quad
\overline{b}_{int}(t^{**})=b^*,
\]
\[
\frac{d}{dt}\tilde a_i(t)=\frac{d}{dt}\tilde b_i(t)=\tilde u_i,\quad 
\tilde a_i(t^{**})=x^*,\quad
\tilde b_i(t^{**})=b^*.
\] 

\begin{Lemma}\label{CNS elargissement}
The previous dynamic satisfies \eqref{CPGD1}-\eqref{CPGD_cons1} if and only if 
\[p_{int}=u_i+p_i-u_{int},\quad
(\overline{u}_{int},\overline{p}_{int})=\left(\frac{u_{int}}{2},\alpha(u_{i}+p_{i})-\frac{u_{int}}{2}\right)
,\quad \tilde p_i=\alpha(u_{i}+p_{i})-\tilde{u}_i,\]
with 
\[0\leq u_{int}\leq u_i+p_i,\quad 0\leq \tilde u_i\leq \alpha(u_i+p_i).\]
\end{Lemma} 

\noindent\underline{Proof}:
We have, for any continuous function $S$ and any function
$\varphi\in\mathcal{D}(\Omega)$,

$$\begin{array}{llll}
& &\hspace{-15mm}\langle\,\,\del_t(nS(u,p,I_{\alpha})+\del_x(nuS(u,p,I_{\alpha})),\varphi\,\,\rangle\\
 &= &-\ds\int_0^{+\infty}\!\!\!\int_{\r}n(t,x)S(u(t,x),p(t,x),I_{\alpha}(x))(\del_t\varphi+u\del_x\varphi)dxdt\\
&=  &-S(u_i,p_i,1)J(0,t^{*},a_i,b_i,u_i)-S(u_{int},p_{int},1)J(t^*,t^{**},a_{int},x^*,u_{int})\\\\
&& -2 S(\overline u_{int},\overline p_{int},1/\alpha)J(t^{*},t^{**},x^*,\overline b_{int},\overline u_{int})\\\\
&&+2 S(\tilde u_{i},\tilde p_{i},1/\alpha)J(t^{**},\infty,\tilde a_i,\tilde b_i,\tilde u_i)\\
&=&-S(u_i,p_i,1) \ds \int_{a^*}^{x^*}\varphi(t^*,x)dx\\
& & -\ds S(u_{int},p_{int},1)\left(-\int_{a^*}^{x^*}\varphi(t^*,x)dx+u_{int}\int_{t^*}^{t^{**}}\varphi(t,x^*)dt\right)\\
&& \ds-2 S(\overline u_{int},\overline p_{int},1/\alpha)\left(\int_{x^*}^{b^*}
\varphi(t^{**},x)dx-\overline u_{int}\int_{t^*}^{t^{**}}\varphi(t,x^*)dt\right)\\
&&+2\ds S(\tilde u_{i},\tilde p_{i},1/\alpha)\int_{x^*}^{b^*}\varphi(t^{**},x)dx,
\end{array}$$
thus we have, in $\mathcal{D}'(\Omega)$,
$$\begin{array}[]{llll}
& & \hspace{-15mm}\del_t(nS(u,p,I_{\alpha})+\del_x(nuS(u,p,I_{\alpha}))\\\\ 
& = &  \left(S(u_{int},p_{int},1)-S(u_i,p_i,1)\right)\1_{[a^*,x^*]}
(x)\delta(t-t^*)\\\\
&  & + 2\left(
S(\tilde u_{i},\tilde p_{i},1/\alpha)- S(\overline u_{int},\overline
p_{int},1/\alpha)\right)\1_{[x^*,b^*]}
(x)\delta(t-t^{**})\\\\
&  & + \left(2\overline u_{int} S(\overline u_{int},\overline
p_{int},1/\alpha)-u_{int}S(u_{int},p_{int},1)\right)\1_{[t^*,t^{**}]}
(t)\delta(x-x^*).
\end{array}$$
For $S(u,p,I_{\alpha})=1$, we get
$$\partial_t n +\partial_x (n u)=\left(2\overline
u_{int}-u_{int}\right)\1_{[t^*,t^{**}]}
(t)\delta(x-x^*).$$
For $S(u,p,I_{\alpha})=(u+p)I_{\alpha}$, we get
$$\begin{array}{lll}
& & \hspace{-15mm}\partial_t (n(u+p)I_{\alpha}) +\partial_x (n u(u+p)I_{\alpha})\\\\
& = &  
\left(u_{int}+p_{int}-u_i-p_i\right)\1_{[a^*,x^*]}
(x)\delta(t-t^*)\\\\
&  & + \frac{2}{\alpha}\left(
\tilde u_{i}+\tilde p_{i}-\overline u_{int}-\overline
p_{int}\right)\1_{[x^*,b^*]}
(x)\delta(t-t^{**})\\\\
&  & + \left(2\overline u_{int}\left(\frac{\overline u_{int}+\overline
p_{int}}{\alpha}\right)-u_{int}(u_{int}+p_{int})\right)\1_{[t^*,t^{**}]}
(t)\delta(x-x^*).
\end{array}$$
Therefore, such a function $(n,u,p)$ is a solution of \eqref{CPGD1}-\eqref{CPGD2} if and only if
\[
\left\{
\begin{array}[]{lll}
2\overline u_{int} & = & u_{int}\\
u_{int}+p_{int} & = & u_i+p_i\\
\tilde u_{i}+\tilde p_{i} & = & \overline u_{int}+\overline
p_{int}\\
\overline u_{int}+\overline p_{int} & = & \alpha(u_{int}+p_{int})
\end{array}\right.\quad\Longleftrightarrow\quad
\left\{
\begin{array}[]{lll}
2\overline u_{int} & = & u_{int}\\
u_{int}+p_{int} & = & u_i+p_i\\
\overline u_{int}+\overline p_{int} & = & \alpha(u_{i}+p_{i})\\
\tilde u_{i}+\tilde p_{i} & = & \alpha(u_{i}+p_{i})
\end{array}\right.,
\]
and we conclude as in lemma \ref{CNS retrecissement}.
\cqfd

Now, if a single block $(n^{*},u_{i},0)$ undergoes a enlargement of the road,
we have 
\[p_{i}=0,\qquad \tilde{p}_{i}=0\]
in the relations of lemma \ref{CNS elargissement}, which leads to 
\[p_{int}=u_i-u_{int},\quad
(\overline{u}_{int},\overline{p}_{int})=\left(\frac{u_{int}}{2},\alpha u_{i}-\frac{u_{int}}{2}\right)
,\quad \tilde u_{i}=\alpha u_{i},\]
with 
\[0\leq u_{int}\leq u_i.\]
In this case, we can't impose $\overline{p}_{int}=0$ since it would imply 
$u_{int}=2\alpha u_{i}$ and then $p_{int} <0$ which is impossible.
Then, we use the second criteria of section \ref{uniqueness dynamics} which is the minimization 
of $p_{int}$ in this case. Here $p_{int}=0$ is possible, that is the choice
$u_{int}=u_{i}$, and the dynamics for a enlargement is the following:

\begin{center}
\includegraphics[width=0.6\textwidth]{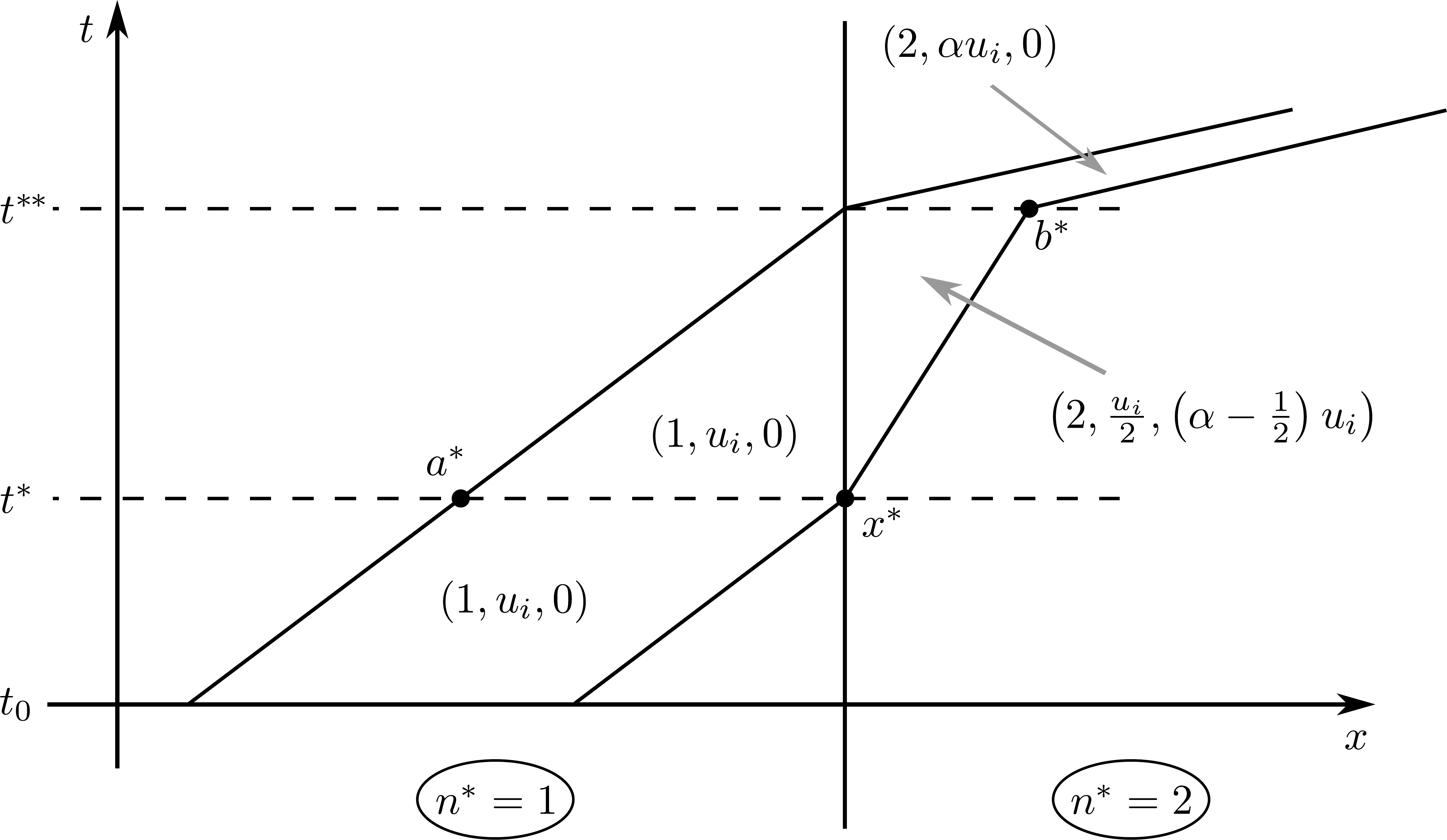}
\end{center} 
This situation is the only case where (\ref{cont+}) can't be additionally asked.
The physical explanation is that the increasing of the speed of the car going from $u_i$ to $\alpha u_i$
is not instantaneous and has to be in two steps. Thus, in the intermediate state, the car is not yet at its preferred velocity and there is 
still a $p$ term.

\begin{Rk}
We notice that the dynamics for enlargement is exactly the reverse process of the narrowing. Moreover, the ML-CPGD model allows cars to accelerate, which was not the case in \cite{berthelin_degond_delitala_rascle} and \cite{berthelin_moutari} (where a maximum principle held for the velocity $u$).
\end{Rk}

\subsection{Compatibility of the dynamics}

Since the previous dynamics are not instantaneous, they can interact before they are completed.
In this subsection, we present the various compatibilities between these dynamics.
Note that it is not just a superposition of various cases.
In order to simplify the presentation, we only show figures that describe the various interactions.

\subsubsection{A train of blocks undergoes an narrowing}

\begin{center}
	\includegraphics[width=0.6\textwidth]{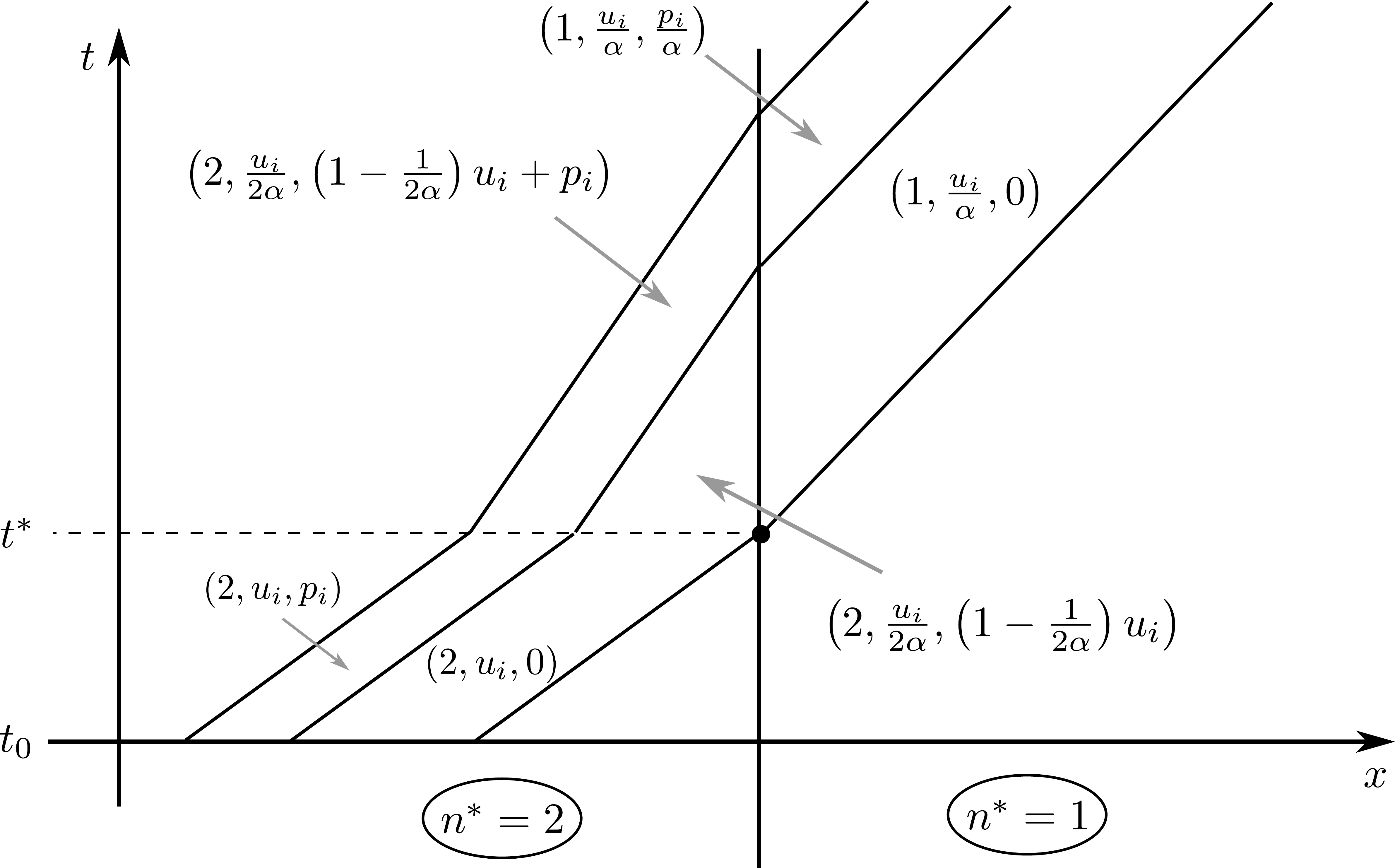}
\end{center}

\subsubsection{Two blocks collide just before the road narrows}

Case with $u_{i-1}>u_{i}$:
\begin{center}
\includegraphics[width=0.6\textwidth]{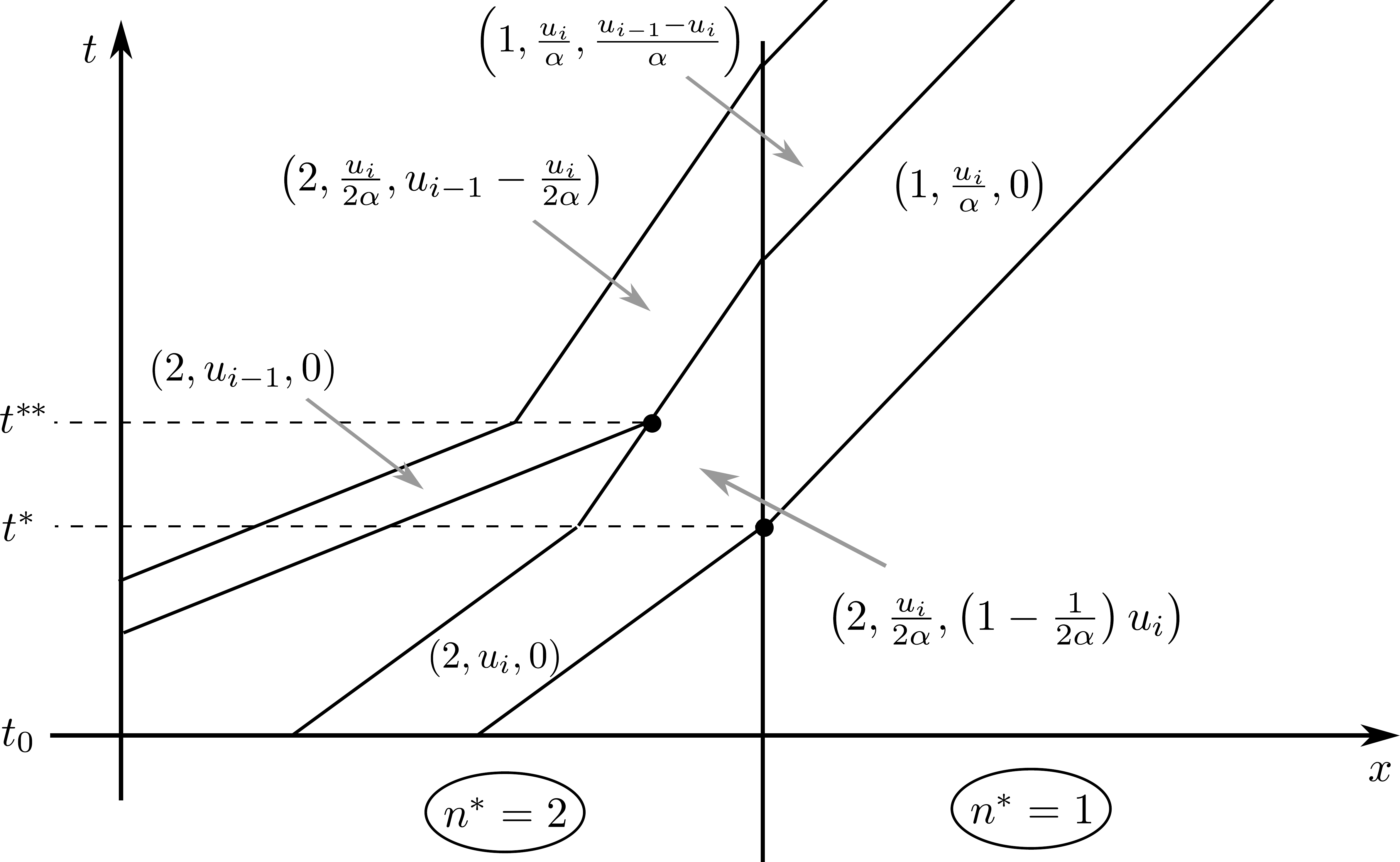}
\end{center} 
Case with $\frac{u_{i}}{2\alpha}<u_{i-1}\leq u_{i}$:
\begin{center}
\includegraphics[width=0.6\textwidth]{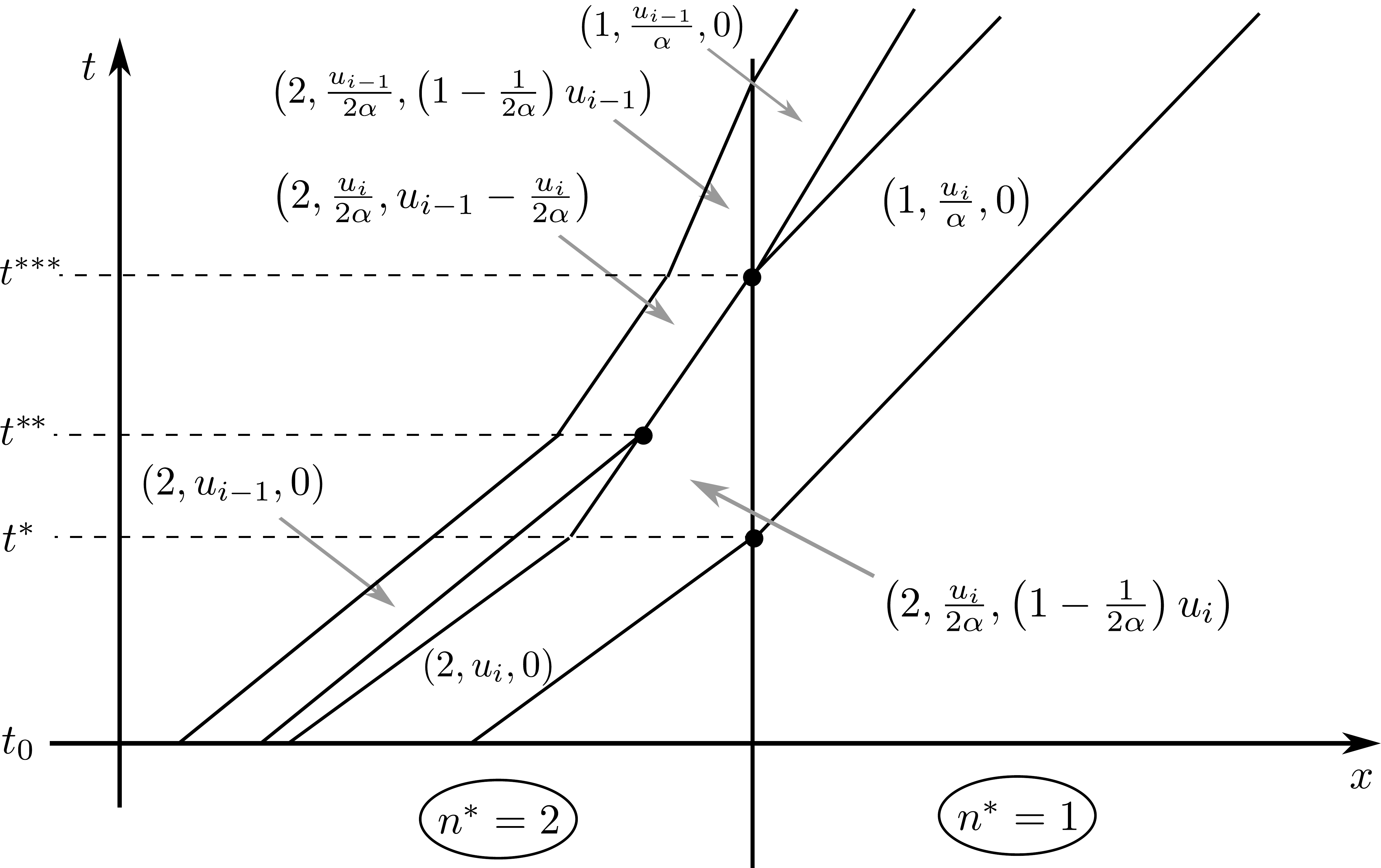}
\end{center} 
Note that in this case, there is
creation of a void area in a jam  due to the fact that 
the acceleration of the leading car is not necessarily followed if there is not a sufficient reserve of speed.
It represents also an approach to model some kind of stop and go waves which is new in such model. 
It will be also the case in the following situation.

\subsubsection{A train of blocks undergoes an enlargement}

\begin{center}
	\includegraphics[width=0.6\textwidth]{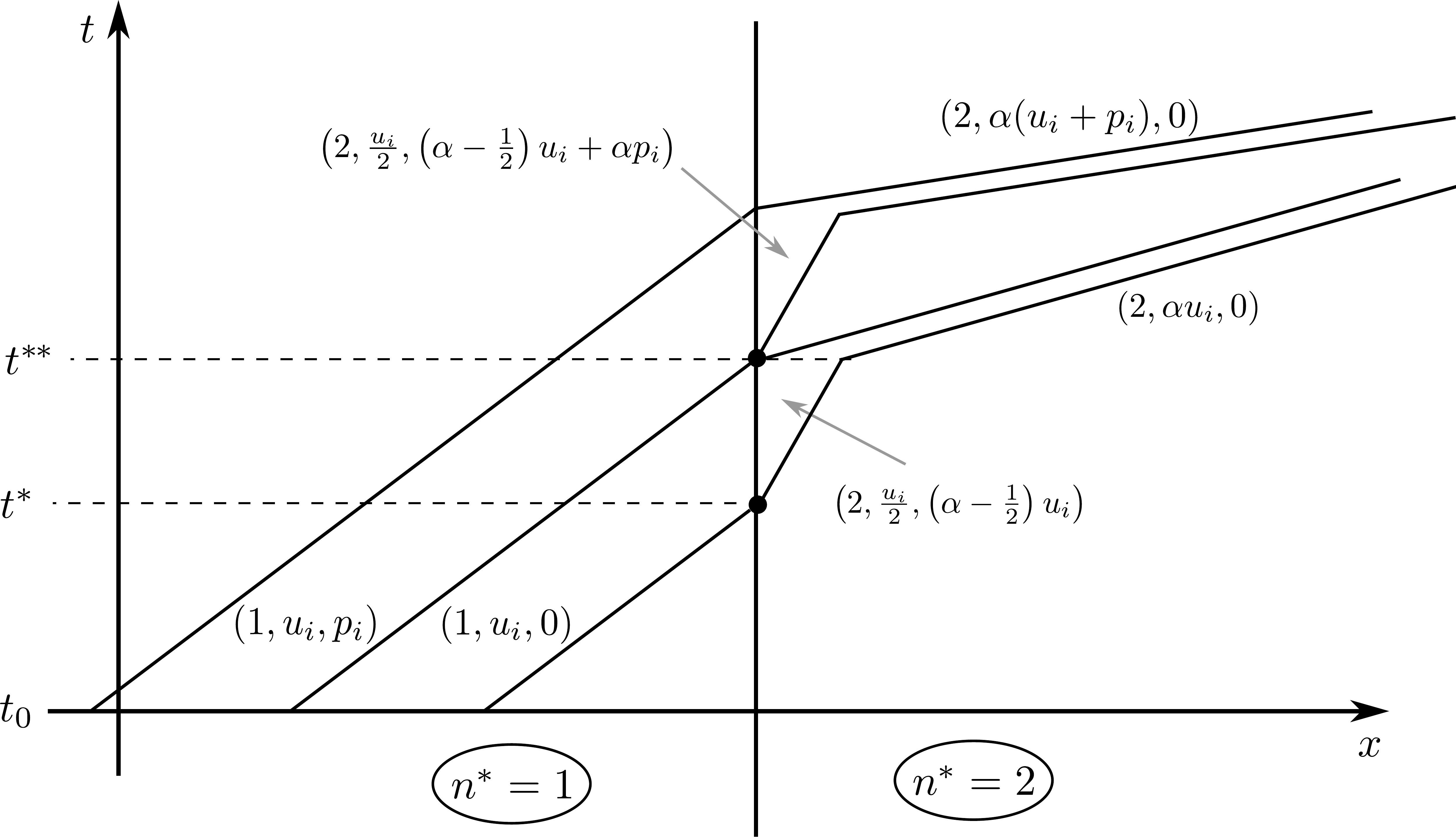}
\end{center}

\subsubsection{Two blocks collide just after the road widens} 

\begin{center}
\includegraphics[width=0.6\textwidth]{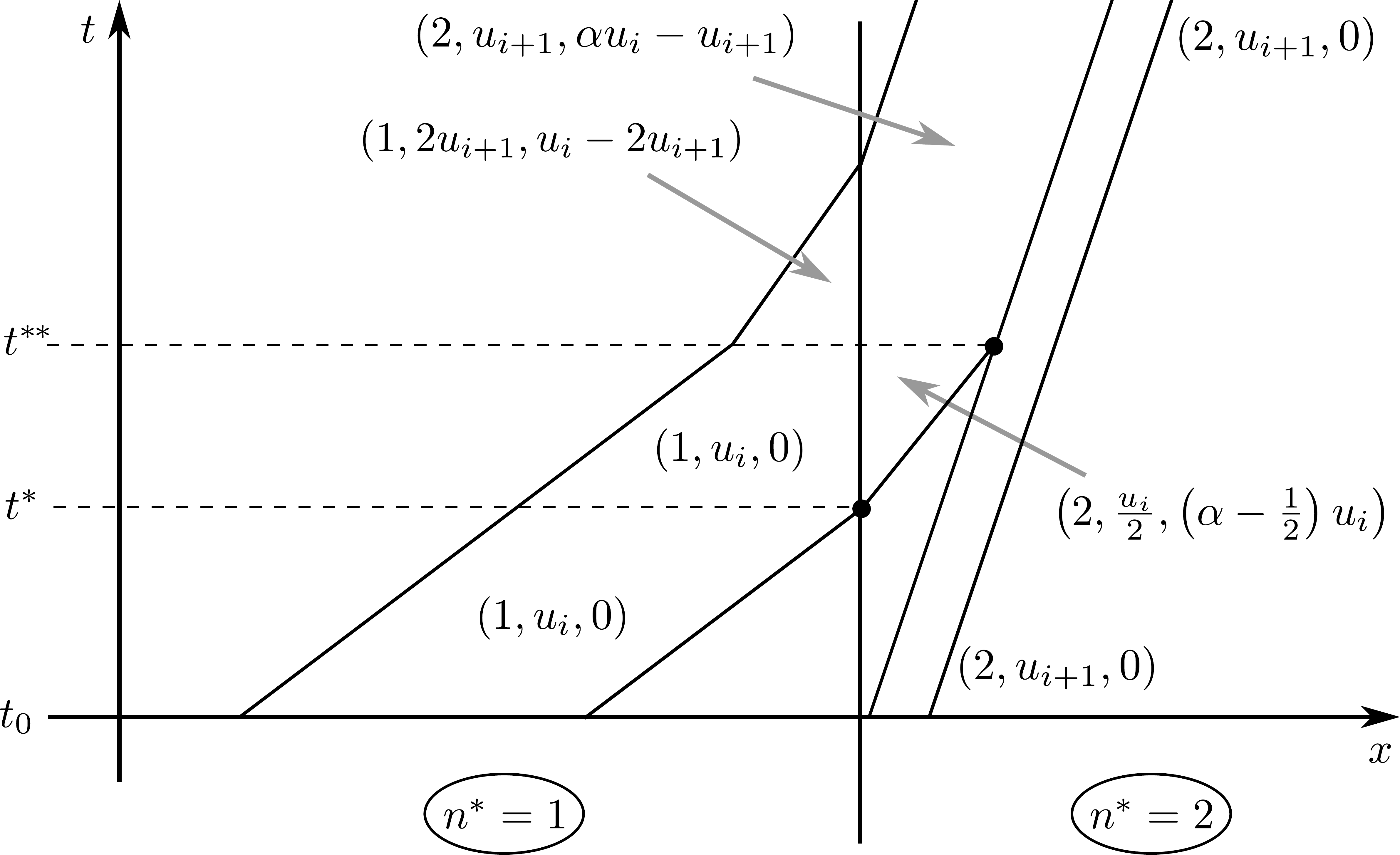}
\end{center} 

\noindent Here, we have $\frac{u_{i}}{2}>u_{i+1}$, thus 
\[\alpha u_{i}-u_{i+1}\geq u_{i}-2u_{i+1}>0.\]

\subsubsection{The road follows $1\to 2\to 1$ faster than the block}

\begin{center}
\includegraphics[width=0.6\textwidth]{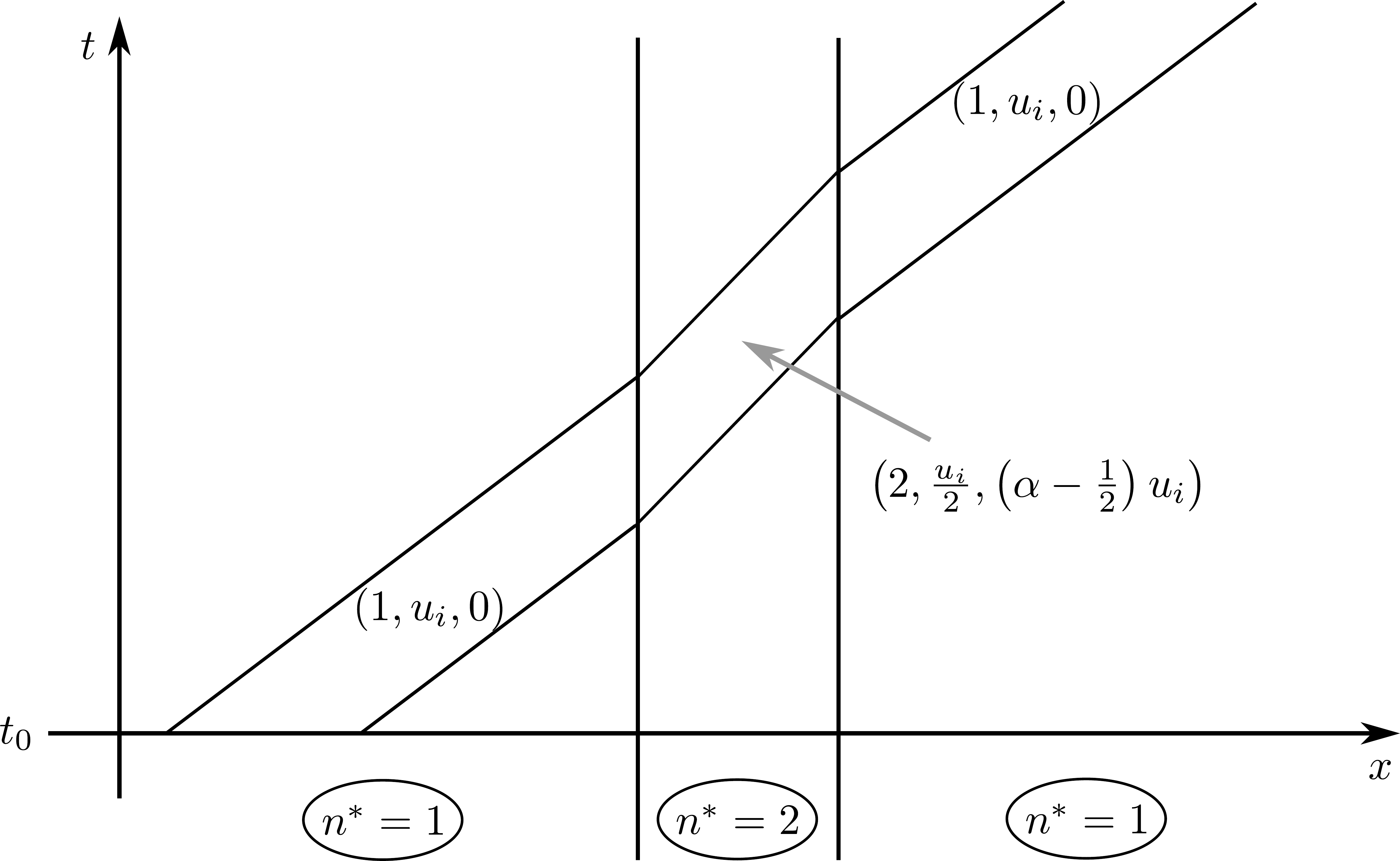}
\end{center} 

\subsubsection{The road follows $2\to 1\to 2$ faster than the block}

\begin{center}
\includegraphics[width=0.6\textwidth]{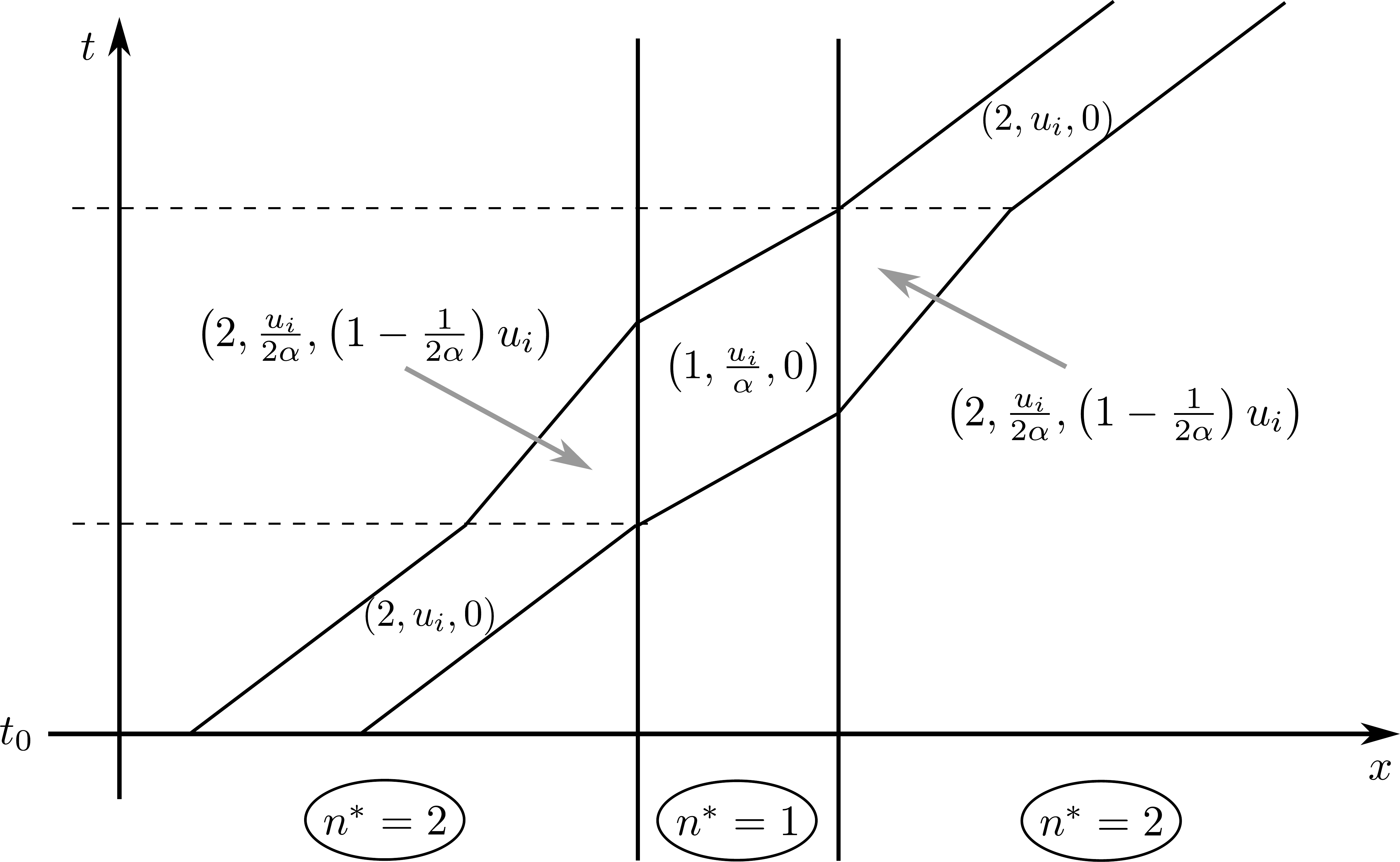}
\end{center} 
\subsection{Block solutions and bounds}

Using the above sections, we are able to state some results on the dynamics of blocks.

\begin{Rk}\label{prolongement_u_p}
The velocity $u$ is assumed to be extended linearly
in the vacuum (areas such that $n=0$) between two successive blocks. Moreover,
we assume that $u$ is constant at $\pm \infty$. But concerning $p$, the
constraint $(n^*-n)p=0$ implies that $p=0$ in the vacuum, and at $\pm\infty$.
Thus, the computations of total variation in $x$ of $u$ and $p$ are different.
\end{Rk}
The previous computations show the following results:
\begin{Th} With the various above dynamics, the quantities $n(t,x)$, $u(t,x)$ and
$p(t,x)$ defined by (\ref{defn})-(\ref{defnp}) and Remark \ref{prolongement_u_p} are solutions to
\eqref{CPGD1}, \eqref{CPGD2}, \eqref{CPGD_cons1}.
\end{Th}
We can also establish some bounds on these solutions:
\begin{Prop} We still denote by $n(t,x)$, $u(t,x)$ and
$p(t,x)$ the functions of (\ref{defn})-(\ref{defnp}) and Remark \ref{prolongement_u_p}.
These functions satisfy the maximum principle 
\begin{equation}\label{principe_max_u}
\ds0\leq u(t,x)\leq 2\alpha\,  \left(esssup_{y}
u^0(y)+esssup_{y}p^0(y)\right),
\end{equation}
\begin{equation}\label{principe_max_p}
\ds0\leq p(t,x)\leq 2\alpha\,  \left(esssup_{y}
u^0(y)+esssup_{y}p^0(y)\right).
\end{equation}
If we assume furthermore that the initial data in the blocks $u_i^0$ and $p_i^0$ are
$BV$ functions, then we have, for all $t\in [0,T]$,
\begin{equation}\label{BV_u}
TV_K(u(t,.))\leq 4\alpha M
\left(TV_{\tilde{K}}(u^0)+TV_{\tilde{K}}(p^0)+\|u^0\|_{L^{\infty}}\right), 
\end{equation}
\begin{equation}\label{BV_p}
 TV_K(p(t,.))\leq 4\alpha
M\left(TV_{\tilde{K}}(u^0)+TV_{\tilde{K}}(p^0)+\|u^0\|_{L^{\infty}}\right), 
\end{equation}
for any compact $K=[a,b]$ and with
\[\tilde{K}=[a-t\,(esssup_{y}u^0),b-t\,(essinf_{y}u^0)],\] where $TV_K$ (resp.
$TV_{\tilde{K}}$) denotes the total variation on the set $K$ (resp.
$\tilde{K}$), and $M$ is the number of road transitions (supposed to be
finite). 
\end{Prop}

\begin{Rk}
The estimate \eqref{principe_max_u} reflects the fact that cars can accelerate in the present model. In \cite{berthelin_degond_delitala_rascle} and \cite{berthelin_moutari}, we simply have the estimate 
\[essinf_{y}u^{0}(y)\leq u(t,x)\leq esssup_{y}u^{0}(y).\]
\end{Rk}

\noindent\underline{Proof}:
We treat some examples which represent the critical cases.
In these cases, we compute the total variation on $\r$ to simplify the presentation.

\begin{itemize}
 \item[$\bullet$]\underline{Case of collisions without change of width:}

We obtain the bounds corresponding to the classical CPGD model (like in \cite{berthelin_moutari}).
 
We assume the following dynamics: 
at time $t=0$, there are $N$ blocks (denoted by $B_1,\dots, B_N$) with velocities $u_1^0>u_2^0>\dots >u_N^0$ 
(which is the case with the most collisions) and pressures $p_1^0,\dots,p_N^0\geq 0$, thus
$$TV(u^0)=\sum_{i=1}^{N-1}|u_i^0-u_{i+1}^0|,\qquad TV(p^0)=2\sum_{i=1}^N p_i^0.$$

 Let $t>0$ such that in the time interval $[0,t]$, the $j$ first blocks $B_{1},\dots, B_{j}$ collide
 successively at $t_1<\dots<t_{j-1}\leq t$ for instance (i.e. $B_{i}$ collide with $B_{i+1}$ at the time $t_{i}$ for all 
 $1\leq i\leq j-1$) 
 and then the $q-j+1$ following blocks $B_{j},\dots,B_{q}$ 
 collide at the same time $t_j$, with $t_{j-1}<t_{j}\leq t$. At the time $t$, the last 
 $N-q+1$ blocks $B_{q},\dots, B_{N}$ have not collided yet.

We have the relations

$$\forall k\in\{1,\dots,j\},\quad \forall i\in\{1,\dots,N\},\qquad p_{i}^{t_{k}}+u_{i}^{t_{k}}=p_{i}^{0}+u_{i}^{0},$$
$$\forall i\in\{1,\dots,q\},\qquad u_{i}^{t_{j}}=u_{q}^{0},\qquad p_{i}^{t_{j}}=u_{i}^{0}+p_{i}^{0}-u_{q}^{0},$$
$$\forall i\in\{q+1,\dots,N\},\qquad u_{i}^{t_{j}}=u_{i}^{0},\qquad p_{i}^{t_{j}}=p_{i}^{0}.$$

Then we get

$$\begin{array}{llllll}
TV(u(t,.)) & = & |u_1^{t_j}-u_2^{t_j}|+\dots+|u_{j-1}^{t_j}-u_{j}^{t_j}|\\\\ 
& & + \quad |u_j^{t_j}-u_{j+1}^{t_j}|+\dots+ \quad|u_{q-1}^{t_j}-u_q^{t_j}|\\\\
&  & + \quad |u_{q}^{t_j}-u_{q+1}^{t_j}|+\dots+|u_{N-1}^{t_j}-u_N^{t_j}|\\\\
& = & \,\,\,\,|u_q^{0}-u_q^{0}|+\dots+|u_{q}^{0}-u_{q}^{0}|\\\\
& & +\quad |u_q^{0}-u_{q+1}^{0}|+\dots+|u_{N-1}^{0}-u_N^{0}|\\\\
& \leq & TV(u^0).
  \end{array}$$

Now, for $p$, we have

$$\begin{array}{llllll}
TV(p(t,.)) & = & p_1^{t_j}+|p_1^{t_j}-p_2^{t_j}|+\dots+|p_{j-1}^{t_j}-p_{j}^{t_j}|\\\\ 
& & +\quad |p_j^{t_j}-p_{j+1}^{t_j}|+\dots+|p_{q-1}^{t_j}-p_{q}^{t_j}|\\\\
 &  & +\quad p_{q}^{t_{j}}\quad +\quad 2\left(p_{q+1}^{t_{j}}+\dots +p_{N}^{t_{j}}\right)\\\\
& = & u_{1}^{0}+p_{1}^{0}-u_{q}^{0}+|u_{1}^{0}+p_{1}^{0}-u_{2}^{0}-p_{2}^{0}|\\\\ 
& & +\dots+|u_{j-1}^{0}+p_{j-1}^{0}-u_{j}^{0}-p_{j}^{0}|\quad\\\\
 & & +\quad |u_{j}^{0}+p_{j}^{0}-u_{j+1}^{0}-p_{j+1}^{0}|\\\\ 
 & & +\dots+|u_{q-1}^{0}+p_{q-1}^{0}-u_{q}^{0}-p_{q}^{0}|\\\\
 &  & +\quad p_{q}^{0}\quad +\quad 2\left(p_{q+1}^{0}+\dots +p_{N}^{0}\right)\\\\
& \leq & u_{1}^{0}-u_{q}^{0}+|u_{1}^{0}-u_{2}^{0}|+\dots+|u_{q-1}^{0}-u_{q}^{0}|\\\\
& & +\quad  2\left(p_{1}^{0}+p_{2}^{0}+\dots +p_{N}^{0}\right)\\\\
& \leq & 2TV(u^0)+TV(p^0).
 \end{array}$$

\item[$\bullet$] \underline{Case of enlargement of the road without collision:} 

We assume the following dynamics: at time $t=0$, we consider two blocks $B_1=(u_1^0,p_1^0)$ and $B_2=(u_2^0,p_2^0)$, in a section of road where $n^*=1$.

We have
$$TV(u^0)=|u_2^0-u_1^0|,\qquad TV(p^0)=2(p_2^0+p_1^0).$$

At time $t_1>0$, the block $B_1$ reach the two-lane section, and undergoes the change of width during the time interval
$[t_1,t_2]$. Then, later in the interval $[t_3,t_4]$ (with $t_3>t_2$) the block $B_2$ enter in the two-lane section.

For all $t\in ]t_1,t_2[$, we have (with the notations of section \ref{elargissement})
$$\begin{array}{lllll}
TV(u(t,.)) & = & |u_2^0-u_{1,int}|+|u_{1,int}-\overline{u_{1,int}}|\\\\ 
& = & |u_2^0-u_{1,int}|+\ds \frac{u_{1,int}}{2}\\\\
& \leq & |u_2^0-u_1^0|+|u_1^0-u_{1,int}|+\ds \frac{u_{1,int}}{2}.
  \end{array}$$

Since $u_{1,int}\leq u_1^0$ we obtain
$$\begin{array}{lllllll}
TV(u(t,.)) & \leq & |u_2^0-u_1^0|+u_1^0-\ds \frac{u_{1,int}}{2} & \leq & TV(u^0)+\|u^0\|_{L^{\infty}}.
  \end{array}$$

Moreover
$$\begin{array}{lllll}
TV(p(t,.)) & = & 2p_2^0+p_{1,int}+|p_{1,int}-\overline{p_{1,int}}|+\overline{p_{1,int}},
  \end{array}$$

but 
\[\overline{p_{1,int}}-p_{1,int}=(\alpha-1)(u_1^0+p_1^0)+\ds \frac{u_{1,int}}{2}\geq 0,\]
thus
\[\begin{array}{lllll}
TV(p(t,.)) & = & 2(p_2^0+\overline{p_{1,int}}) & = & 2(p_2^0+\alpha(u_1^0+p_1^0)-\ds \frac{u_{1,int}}{2}),
\end{array}
\]

and we deduce
$$\begin{array}{lllllll}
TV(p(t,.)) & \leq & 2\alpha(p_1^0+p_2^0)+2\alpha u_1^0 & \leq & \alpha TV(p^0)+2\alpha \|u^0\|_{L^{\infty}}.
\end{array}$$

For all $t\in ]t_2,t_3[$, we have 
$$\begin{array}{lllllll}
TV(u(t,.)) & = & |u_2^0-\alpha u_1^0| & \leq & |u_2^0-u_1^0|+|(1-\alpha)u_1^0|\\\\
& &  & \leq & TV(u^0)+(\alpha-1)\|u^0\|_{L^{\infty}},
  \end{array}$$

and
\[
TV(p(t,.))  =  2p_2^0+2\alpha p_1^0  \leq  \alpha TV(p^0).
\]

For all $t\in ]t_3,t_4[$, we have 
\[\begin{array}{lllllll}
TV(u(t,.)) & = & |u_{2,int}-\overline{u_{2,int}}|+|\overline{u_{2,int}}-\alpha u_1^0|\\\\
& = &\ds \frac{u_{2,int}}{2}+|\frac{u_{2,int}}{2}-\alpha u_1^0|\\\\
& \leq &\ds \frac{u_{2,int}}{2}+|\frac{u_{2,int}}{2}-\alpha u_2^0|+|\alpha u_2^0-\alpha u_1^0|.
\end{array}\]

Since $u_{2,int}\leq u_2^0$, we obtain

$$\begin{array}{lllllll}
TV(u(t,.)) & \leq &\ds \frac{u_{2,int}}{2}-\frac{u_{2,int}}{2}+\alpha u_2^0+|\alpha u_2^0-\alpha u_1^0|\\\\ 
& \leq & \alpha TV(u^0)+\alpha \|u^0\|_{L^{\infty}}.
\end{array}$$

Moreover

$$\begin{array}{lllllll}
TV(p(t,.)) & = & p_{2,int}+|p_{2,int}-\overline{p_{2,int}}|+\overline{p_{2,int}}+2\alpha p_1^0\\\\ 
& = & 2(\overline{p_{2,int}}+\alpha p_1^0)\\\\
& = & 2(\alpha(u_2^0+p_2^0)-\ds\frac{u_{2,int}}{2}+\alpha p_1^0),
  \end{array}$$

and we deduce
$$\begin{array}{lllllll}
TV(p(t,.)) & \leq & 2\alpha(p_1^0+p_2^0)+2\alpha u_2^0 & \leq & \alpha TV(p^0)+2\alpha \|u^0\|_{L^{\infty}}.
\end{array}$$

At least, for $t>t_4$, we have
$$\begin{array}{lllllll}
TV(u(t,.)) & = & |\alpha u_2^0-\alpha u_1^0| & = & \alpha TV(u^0),
  \end{array}$$

and 
$$\begin{array}{lllllll}
TV(p(t,.)) & = & 2\alpha p_2^0+2\alpha p_1^0 & = & \alpha TV(p^0).
  \end{array}$$

Finally, the bound is
\[TV(u(t,.))\leq \alpha(TV(u^0)+\|u^0\|_{L^{\infty}}),\]
\[TV(p(t,.))\leq \alpha(TV(u^0)+2\|u^0\|_{L^{\infty}}),\] 
for all $t>0$.

In the general case (if we follow $N$ blocks along the time), we shall obtain the same bound because only one block at a time undergoes every change $n^*=1\rightarrow 2$.

But it is possible that many blocks undergo this enlargement together at different places. That is why the general estimate is the following:

$$\begin{array}{llll}
TV(u(t,.)) & \leq & \alpha(TV(u^0)+M\|u^0\|_{L^{\infty}}),
  \end{array}$$
$$\begin{array}{llll}
TV(p(t,.)) & \leq & \alpha (TV(p^0)+2M\|u^0\|_{L^{\infty}}),
  \end{array}$$

where $M$ is the number of lane transitions.

 \item[$\bullet$]\underline{Case of narrowing of the road without collision:}

The computations are similar to the previous case. With the notations of section \ref{retrecissement}, we have:

For all $t\in ]t_1,t_2[$, 
$$\begin{array}{lllll}
TV(u(t,.)) & = & |u_2^0-u_{1,int}|+|u_{1,int}-\overline{u_{1,int}}|\\\\ 
& = & |u_2^0-u_{1,int}|+u_{1,int}\\\\
& \leq & |u_2^0-u_1^0|+|u_1^0-u_{1,int}|+u_{1,int}.
  \end{array}$$

Since $u_{1,int}\leq u_1^0$ we obtain
$$\begin{array}{lllllll}
TV(u(t,.)) & \leq & |u_2^0-u_1^0|+u_1^0 & \leq & TV(u^0)+\|u^0\|_{L^{\infty}}.
  \end{array}$$

Moreover,
$$\begin{array}{lllll}
TV(p(t,.)) & = & 2p_2^0+p_{1,int}+|p_{1,int}-\overline{p_{1,int}}|+\overline{p_{1,int}},
  \end{array}$$

but this time, $\overline{p_{1,int}}\leq p_{1,int}$, thus
$$\begin{array}{lllll}
TV(p(t,.)) & = & 2(p_2^0+p_{1,int}) & = & 2(p_2^0+u_1^0+p_1^0-u_{1,int}),
  \end{array}$$

and we deduce
$$\begin{array}{lllllll}
TV(p(t,.)) & \leq & 2(p_1^0+p_2^0)+2 u_1^0 & \leq & TV(p^0)+2 \|u^0\|_{L^{\infty}}.
\end{array}$$

For all $t\in ]t_2,t_3[$, we have 
$$\begin{array}{lllllll}
TV(u(t,.)) & = & |u_2^0-\ds\frac{1}{\alpha} u_1^0|\quad \leq \quad|u_2^0-u_1^0|+|(1-\ds\frac{1}{\alpha})u_1^0|\\\\ 
& \leq & TV(u^0)+(1-\ds\frac{1}{\alpha})\|u^0\|_{L^{\infty}},
  \end{array}$$

and
$$\begin{array}{lllllll}
TV(p(t,.)) & = & 2p_2^0+\ds\frac{2}{\alpha} p_1^0 & \leq &  TV(p^0).
  \end{array}$$

For all $t\in ]t_3,t_4[$, we have 
$$\begin{array}{lllllll}
TV(u(t,.)) & = & |u_{2,int}-\overline{u_{2,int}}|+|\overline{u_{2,int}}-\ds\frac{1}{\alpha} u_1^0|\\\\ 
& =& u_{2,int}+|2u_{2,int}-\ds\frac{1}{\alpha} u_1^0|\\\\
& \leq & u_{2,int}+|2u_{2,int}-\ds\frac{1}{\alpha} u_2^0|+\ds|\frac{1}{\alpha} u_2^0-\ds\frac{1}{\alpha} u_1^0|.
  \end{array}$$

Since $u_{2,int}\leq u_2^0$, we obtain

$$\begin{array}{lllllll}
TV(u(t,.)) & \leq & \ds u_{2,int}+2u_2^0-2u_{2,int}+(2-\frac{1}{\alpha})u_2^0+\frac{1}{\alpha}|u_2^0-u_1^0|\\\\
& \leq & 4u_2^0 + \ds\frac{1}{\alpha}TV(u^0)\\\\
& \leq & 4\|u^0\|_{L^{\infty}}+\ds\frac{1}{\alpha}TV(u^0).
  \end{array}$$
Moreover
$$\begin{array}{lllllll}
TV(p(t,.)) & = & p_{2,int}+|p_{2,int}-\overline{p_{2,int}}|+\overline{p_{2,int}}+\ds\frac{2}{\alpha} p_1^0\\\\ 
& = & 2(p_{2,int}+\ds\frac{1}{\alpha} p_1^0)\\\\
 & = & 2(u_2^0+p_2^0-u_{2,int}+\ds\frac{1}{\alpha} p_1^0),
  \end{array}$$

and we deduce
$$\begin{array}{lllllll}
TV(p(t,.)) & \leq & 2(p_1^0+p_2^0)+2 u_2^0 & \leq & TV(p^0)+2 \|u^0\|_{L^{\infty}}.
\end{array}$$

At least, for $t>t_4$, we have
$$\begin{array}{lllllll}
TV(u(t,.)) & = &\ds  |\frac{1}{\alpha} u_2^0-\frac{1}{\alpha} u_1^0| & = & \ds\frac{1}{\alpha} TV(u^0),
  \end{array}$$

and 
$$\begin{array}{lllllll}
TV(p(t,.)) & = &\ds \frac{2}{\alpha} p_2^0+\frac{2}{\alpha} p_1^0 & = & \ds\frac{1}{\alpha} TV(p^0).
  \end{array}$$

Finally, the bound is
\[TV(u(t,.))\leq TV(u^0)+4\|u^0\|_{L^{\infty}},\]
\[TV(p(t,.))\leq TV(u^0)+2\|u^0\|_{L^{\infty}},\] 
for all $t>0$.

Now $\|u^0\|_{L^{\infty}}$ can appear on every lane transition and the  estimate is then

$$\begin{array}{llll}
TV(u(t,.)) & \leq & TV(u^0)+4M\|u^0\|_{L^{\infty}},
  \end{array}$$
$$\begin{array}{llll}
TV(p(t,.)) & \leq & TV(p^0)+2M\|u^0\|_{L^{\infty}},
  \end{array}$$

where $M$ is the number of lane transitions.

 \item[$\bullet$]The general situation is a superposition of these cases and it gives the Proposition.\cqfd
\end{itemize}

\section{Existence of weak solutions}

In this section, we prove the existence of weak solutions using previous 
clusters dynamics, an approximation lemma of the initial data by these sticky blocks and a
compactness result.

\subsection{Approximation of the initial data by sticky blocks}

We have the following lemma, which is widely inspired from the ones in \cite{bertheli}
and \cite{berthelin_moutari}, but here, $n^*$ is piecewise constant, and
constant at $\pm\infty$. We can see the Appendix for the proof of this variant.
\noindent

\begin{Lemma}\label{lemme_approximation} Let $n^0\in L^1(\r)$, $u^0,p^0\in
L^{\infty}(\r)\cap BV(\r)$ such
that $0\leq n^0\leq n^*(x)$, $0\leq u^0$, $0\leq p^0$ and $(n^*(x)-n^0)p^0=0$.
Then, there exists a
sequence of block initial data $(n_k^0,u_k^0,p_k^0)_{k\geq 1}$ such that for all $k\in\N^{*}$,

\begin{equation}\label{approx_1}
 \int_{\r}n_k^0(x)dx\leq\int_{\r}n^0(x)dx,
\end{equation}

\begin{equation}\label{approx_2}
essinf\,\,u^0\leq u_k^0\leq esssup\,\,u^0,\quad essinf\,\,p^0\leq p_k^0\leq
esssup\,\, p^0,
\end{equation}

\begin{equation}\label{approx_3}
TV(u_k^0)\leq TV(u^0),\quad TV(p_k^0)\leq TV(p^0),
\end{equation}
and for which the convergences $n_k^0\tw n^0$, $n_k^0u_k^0\tw n^0u^0$ and
$n_k^0p_k^0\tw n^0p^0$ hold in the distribution sense. 
Moreover, the sequence $(n_k^0,u_k^0,p_k^0)$ satisfies the constraint:
$$(n^*(x)-n_k^0)p_k^0=0,\qquad\forall k\geq 1.$$
\end{Lemma}

\subsection{Existence result}

Let us recall the ML-CPGD system:
\begin{eqnarray}
& &       \partial_t n +  \partial_x (n u) = 0 \, ,  \label{CPGD1bis}\\
& &      \partial_t (n(u+p)I_{\alpha})+  \partial_x (n u(u + p )I_{\alpha}) = 0
\, , \label{CPGD2bis}\\
& & 0 \leq n \leq n^*(x) \, , \quad \, u\geq 0 \,, \quad p \geq 0 \, , \quad 
(n^*(x) - n) p = 0 \,.
  \label{CPGD_cons1bis}
\end{eqnarray}

We prove now the existence of weak solutions.
The idea is first to approximate the initial data in the distributional sense by sticky blocks.
These special initial data give a sequence of solutions.
Then we perform a compactness argument on this sequence of solutions.
Finally, we prove that the obtained limit is a solution for the wanted initial data.
The regularity of the solutions are
\begin{equation}\label{regularite n}
n\in L^{\infty}(]0,+\infty[_t,L^{\infty}(\r_x)\cap L^1(\r_x)),
\end{equation}
\begin{equation}\label{regularite u et p}
u,p\in L^{\infty}(]0,+\infty[_t,L^{\infty}(\r_x)).
\end{equation}

\begin{Th}\label{compacite sur le systeme cpgd}

Let $(n^0,u^0,p^0)$ be some initial data such that 
\[n^0\in L^1(\r),\qquad u^0,p^0\in L^{\infty}(\r)\cap BV(\r),\] with
$0\leq u^0$, $0\leq p^0$, $0\leq n^0\leq n^*(x)$ and $(n^*(x)-n^0)p^0=0$.
Then there exists $(n,u,p)$ with regularities \eqref{regularite n},
\eqref{regularite u et p}, solution to the system
$\eqref{CPGD1bis}-\eqref{CPGD_cons1bis}$, with initial data $(n^0,u^0,p^0)$.
The obtained solution also satisfies
\begin{equation}\label{principe_max_u_bis}
\ds0\leq u(t,x)\leq 2\alpha\,  \left(esssup_{y}
u^0(y)+esssup_{y}p^0(y)\right),
\end{equation}
\begin{equation}\label{principe_max_p_bis}
\ds0\leq p(t,x)\leq 2\alpha\,  \left(esssup_{y}
u^0(y)+esssup_{y}p^0(y)\right).
\end{equation}
\end{Th}

\noindent\underline{Proof}:
Let $n_k^0$, $u_k^0$, $p_k^0$ ($k\in\N^*$) be the block
initial data associated respectively to $n^0$, $u^0$, $p^0$ provided by Lemma
\ref{lemme_approximation}. For all $k$, the results of section
\ref{bouchons_collants} allow us to get $(n_k,u_k,p_k)$ solutions of
$\eqref{CPGD1bis}-\eqref{CPGD_cons1bis}$ with initial data
$(n_k^0,u_k^0,p_k^0)$, with regularities \eqref{regularite n},
\eqref{regularite u et p}, and which satisfy the bounds
\begin{equation}\label{principe_max_u_k}
\ds0\leq u_k(t,x)\leq 2\alpha\,  \left(esssup_{y}
u^0_k(y)+esssup_{y}p^0_k(y)\right),
\end{equation}
\begin{equation}\label{principe_max_p_k}
\ds0\leq p_k(t,x)\leq 2\alpha\,  \left(esssup_{y}
u^0_k(y)+esssup_{y}p^0_k(y)\right),
\end{equation}
\begin{equation}\label{BV_u_k}
TV_K(u_k(t,.))\leq 4\alpha M
\left(TV_{\tilde{K}}(u_k^0)+TV_{\tilde{K}}(p_k^0)+\|u^0_k\|_{L^{\infty}}\right),
\end{equation}
\begin{equation}\label{BV_p_k}
 TV_K(p_k(t,.))\leq 4\alpha M
\left(TV_{\tilde{K}}(u_k^0)+TV_{\tilde{K}}(p_k^0)+\|u^0_k\|_{L^{\infty}}\right).
\end{equation}
Since $(n_k)$ is bounded  in $L^{\infty}$, then there exists a
subsequence such that
\begin{equation}
\mbox{$n_k\tw n$ in $L^{\infty}_{w*}(]0,+\infty[\times\r)$}.
\end{equation}
Thanks to \eqref{principe_max_u_k}, \eqref{principe_max_p_k} and the bounds on
$u_k^0,p_k^0$ provided by Lemma \ref{lemme_approximation}, the sequence $(u_k)$
and $(p_k)$ are bounded in $L^{\infty}(]0,+\infty[\times\r)$, then, up to
subsequences, we have
\begin{equation}
\mbox{$u_k\tw u$ in $L^{\infty}_{w*}(]0,+\infty[\times\r)$},
\end{equation}
\begin{equation}
\mbox{$p_k\tw p$ in $L^{\infty}_{w*}(]0,+\infty[\times\r)$}.
\end{equation}
 
Next step is now to prove the passage to the limit in the equation.\\
\noindent First, for the sequence $(n_k)_{k\geq 1},$ we can obtain more
compactness using the following lemma and the estimate:

$\forall T>0,\quad \forall \varphi\in\mathcal{D}(\r_x),\quad \forall
t,s\in[0,T],\quad \forall k\in \N^*,$
\begin{equation}\label{bornes lipschitz en t sur n_k}
\begin{split}
 & \hspace{-20mm}\left|\int_{\r}(n_k(t,x)-n_k(s,x))\varphi(x)dx\right|\\ 
\leq & \quad n^{*}\sup_{k\geq 1}\|u_k^0\|_{L^{\infty}(\r_x)}\left(\int_{\r}
|\del_x\varphi|dx\right)|t-s|,
\end{split}
\end{equation}
which can be obtained by integrating \eqref{CPGD1bis}.
\begin{Lemma}\label{lemme de compacite}
Let $(n_k)_{k\in\N^*}$ be a bounded sequence in $L^{\infty}(]0,T[\times \r)$
which satisfies: 
for all $\varphi\in\mathcal{D}(\r_x)$, the sequence
$\left(\int_{\r}n_k(t,x)\varphi(x)dx\right)_k$ is uniformly Lipschitz
continuous on $[0,T]$, i.e.
\[\exists C_{\varphi}>0,\quad \forall
k\in\N^*,\quad \forall s,t \in[0,T],\]
\[\left|\int_{\r}(n_k(t,x)-n_k(s,x))\varphi(x)dx\right|\leq C_{\varphi}|t-s|.\]
Then, up to a subsequence, it exists $n\in
L^{\infty}(]0,T[\times\r)$ such that $n_k\rightarrow n$ in
$C([0,T],L^{\infty}_{w*}(\r_x))$, i.e.
$$\forall \Gamma\in L^1(\r_x),\quad
\sup_{t\in[0,T]}\left|\int_{\r}(n_k(t,x)-n(t,
x))\Gamma(x)dx\right|\underset{k}{\longrightarrow} 0.$$
\end{Lemma}

\noindent\underline{Proof}:
It is a classical argument of equicontinuity. We can see Appendix for the details. 	
\cqfd
\noindent\underline{Following of the proof of Theorem \ref{compacite sur le systeme
cpgd}}: According to \eqref{CPGD_cons1bis} and \eqref{bornes lipschitz en t sur
n_k}, the lemma \ref{lemme de compacite} applies to the sequence $(n_k)_{k\geq
1}$, and thus 
\begin{equation}\label{limite_n_k}
\mbox{$n_k\rightarrow n$ in $C([0,T],L^{\infty}_{w*}(\r_x))$, for all
$T>0$.}
\end{equation}
As the same, we obtain (integrating
\eqref{CPGD2bis}) an estimate similar to \eqref{bornes lipschitz en t sur
n_k} for the sequence $(n_k(u_k+p_k)I_{\alpha})_{k\geq 1}$, thus it exists $q\in
L^{\infty}(]0,+\infty[\times \r)$ such that 
\begin{equation}\label{limite_n_k(u_k+p_k)I_alpha}
\mbox{$n_k(u_k+p_k)I_{\alpha}\rightarrow
q$ in
$C([0,T],L^{\infty}_{w*}(\r_x))$, for all $T>0$.}
\end{equation}
Now, the key point of the proof is passing to the limit in the products and is
treated by the following technical lemma:

\begin{Lemma}\label{compacite par compensation Florent}
Let us assume that $(\gamma_k)_{k\in\N}$ is a bounded sequence in
$L^{\infty}(]0,T[\times\r)$ that tends to $\gamma$ in
$L^{\infty}_{w*}(]0,T[\times \r),$ and satisfies for any
$\Gamma\in\mathcal{D}(\r_x),$
\begin{equation}\label{condition 1 cc}
\int_{\r}(\gamma_k-\gamma)(t,x)\Gamma(x)dx\underset{k}{\longrightarrow} 0, 
\end{equation}
either i) a.e. $t\in]0,T[$ or ii) in $L^1(]0,T[_t)$. \\
\noindent Let us also assume that $(\omega_k)_{k\in\N}$ is a bounded sequence in
$L^{\infty}(]0,T[\times\r)$ that tends to $\omega$ in
$L^{\infty}_{w*}(]0,T[\times \r),$ and such that for all compact interval
$K=[a,b]$, there exists $C>0$ such that the total variation (in x) of
$\omega_k$ and $\omega$ over $K$ satisfies 
\begin{equation}\label{condition 2 cc}
\forall k\in\N,\quad TV_K(\omega_k(t,.))\leq C,\qquad 
TV_K(\omega(t,.))\leq C.
\end{equation}
Then, $\gamma_k\omega_k\tw \gamma\omega$ in $L^{\infty}_{w*}(]0,T[\times \r),$
as $k\rightarrow +\infty$.
\end{Lemma}

\begin{Rk} This is a result of compensated compactness, which uses the
compactness in $x$ for $(\omega_k)_k$ given by \eqref{condition 2 cc} and the
weak compactness in $t$ for $(\gamma_k)_k$ given by \eqref{condition 1 cc} to
pass to the weak limit in the product $\gamma_k\omega_k$. 
\end{Rk}

\noindent\underline{Proof}: We can refer to \cite{bertheli} for a complete proof, even in the case where
$$
\forall k\in\N,\quad TV_K(\omega_k(t,.))\leq C(1+\frac{1}{t}),\qquad 
TV_K(\omega(t,.))\leq C(1+\frac{1}{t})
,$$
which is more general.\cqfd

\vspace{2mm}

\noindent\underline{End of the proof of Theorem \ref{compacite sur le systeme
cpgd}}:
The convergence \eqref{limite_n_k} allows to apply Lemma \ref{compacite par
compensation Florent} with $\gamma_k=n_k$. Moreover, thanks to
$\eqref{BV_u_k}$ and the $BV$ bounds on $u_k^0$ provided by Lemma
\ref{lemme_approximation}, we can set $\omega_k=u_k$ in Lemma \ref{compacite
par compensation Florent} (in fact, the sequence $u_k(t,.)$ is uniformly
bounded in $BV$ with respect to $t$, and also $u(t,.)$ thanks to the lower
semi-continuity to the $BV$ norm). Thus, we have 
\begin{equation}\label{limite_n_ku_k}
\mbox{$n_ku_k\tw nu$ in
$L^{\infty}_{w*}(]0,+\infty[\times \r)$}.
\end{equation}
The same applies to the sequences $(\gamma_k,\omega_k)=(n_k,p_k)$ and
$(\gamma_k,\omega_k)=(n_k(u_k+p_k)I_{\alpha},u_k)$: we have
\begin{equation}\label{limite_n_kp_k}
\mbox{$n_kp_k\tw np$ in
$L^{\infty}_{w*}(]0,+\infty[\times \r)$},
\end{equation}
\begin{equation}\label{limite_n_ku_k(u_k+p_k)I_alpha}
\mbox{$n_k(u_k+p_k)I_{\alpha}u_k\tw qu$ in
$L^{\infty}_{w*}(]0,+\infty[\times \r)$}.
\end{equation}
Furthermore, we easily have 
\[\mbox{$n_k(u_k+p_k)I_{\alpha}\tw n(u+p)I_{\alpha}$ in $L^{\infty}_{w*}(]0,+\infty[\times \r)$},\] 
thus $q=n(u+p)I_{\alpha}$, and 
\begin{equation}\label{derniere_limite}
\mbox{$n_ku_k(u_k+p_k)I_{\alpha}\tw nu(u+p)I_{\alpha}$ in
$L^{\infty}_{w*}(]0,+\infty[\times \r)$}.
\end{equation}
We deduce that $(n,u,p)$ satisfies $\eqref{CPGD1bis}$, $\eqref{CPGD2bis}$ in
$\mathcal{D}'(]0,+\infty[\times\r)$, and the constraints
$\eqref{CPGD_cons1bis}$.\\
\noindent The last step is to show that $(n^0,p^0,u^0)$ is really the initial data of the
problem, according to the weak formulation:
$$\begin{array}{llll}
\forall\varphi\in C^{\infty}_c([0,+\infty[_t\times\r_x),\\\\
\ds\int_{0}^{\infty}\!\!\!\int_{\r}\left(n\del_t\varphi+nu\del_x
\varphi\right)(t,x)dxdt+\int_{\r}n^0(x)\varphi(0,x)dx=0,\\\\
\ds\int_{0}^{\infty}\!\!\!\int_{\r}\left(n(u+p)I_{\alpha}\del_t\varphi+nu(u+p)I_{
\alpha } \del_x\varphi\right)(t,x)dxdt\\\\
\qquad +\qquad \ds\int_{\r}n^0(x)(u^0(x)+p^0(x))I_{\alpha}(x)\varphi(0,
x)dx=0.
\end{array}$$
It comes easily, because we
have, for all $k\geq 1$:
$$\begin{array}{lll}
\forall
\varphi\in C^{\infty}_c([0,+\infty[_t\times\r_x),\\\\
\ds\int_{0}^{\infty}\!\!\!\int_{\r}\left(n_k\del_t\varphi+n_ku_k\del_x
\varphi\right)(t,x)dxdt+\int_{\r}n_k^0(x)\varphi(0,x)dx=0,\\\\
\ds\int_{0}^{\infty}\!\!\!\int_{\r}\left(n_k(u_k+p_k)I_{\alpha}
\del_t\varphi+n_ku_k(u_k+p_k)I_ {\alpha }\del_x\varphi\right)(t,x)dxdt\\\\
\qquad +\qquad \ds\int_{\r}n_k^0(x)(u_k^0(x)+p_k^0(x))I_{\alpha}
(x)\varphi(0 ,x)dx=0,
\end{array}
$$
and we can pass to the limit when $k\rightarrow +\infty$ because of the
convergences
$n_k^0\tw n^0$, $n_k^0u_k^0\tw n^0u^0$ and $n_k^0p_k^0\tw n^0p^0$ in
$\mathcal{D}'(\r)$, and the convergences
$\eqref{limite_n_ku_k}$, $\eqref{limite_n_ku_k(u_k+p_k)I_alpha}$ and
$\eqref{derniere_limite}$
 in $L^{\infty}_{w*}(]0,+\infty[\times
\r).$
\cqfd

\subsection{Compactness result}

To finalize the paper, we set a compactness result which is contained into the
proof of the previous existence Theorem.

\begin{Th}

Let us consider a sequence of solutions $(n_k,u_k,p_k)$ with regularity \eqref{regularite n},
\eqref{regularite u et p}, satisfying $\eqref{CPGD1bis}-\eqref{CPGD_cons1bis}$, and the following bounds:
$$\forall k\in \N,\quad a.e.\,\, (t,x)\in ]0,+\infty[\times\r,\qquad 0\leq u_k(t,x)\leq C_{\alpha},$$
$$\forall k\in \N,\quad a.e.\,\, (t,x)\in ]0,+\infty[\times\r,\qquad 0\leq p_k(t,x)\leq C_{\alpha},$$
$$\forall K=[a,b]\subset \r,\quad \forall k\in \N,\quad a.e.\,\, t\in ]0,+\infty[,\qquad TV_K(u_k(t,.))\leq C_{\alpha,M,K},$$
$$\forall K=[a,b]\subset \r,\quad \forall k\in \N,\quad a.e.\,\, t\in ]0,+\infty[,\qquad TV_K(p_k(t,.))\leq C_{\alpha,M,K},$$
with $C_{\alpha}$ (resp. $C_{\alpha,M,K}$) some positive constant depending only on $\alpha$ (resp. $\alpha$, $M$ and $K$).\\
\noindent Then, up to a subsequence, $(n_k,u_k,p_k)\tw (n,u,p)$ in $L^{\infty}_{w*}(]0,+\infty[\times\r)$,
where $(n,u,p)$ is a solution to the system
$\eqref{CPGD1bis}-\eqref{CPGD_cons1bis}$.
This solution $(n,u,p)$ also satisfies
$$a.e.\,\, (t,x)\in ]0,+\infty[\times\r,\qquad 0\leq u(t,x)\leq C_{\alpha},$$
$$a.e.\,\, (t,x)\in ]0,+\infty[\times\r,\qquad 0\leq p(t,x)\leq C_{\alpha}.$$
\end{Th}
 
\section{Appendix}

\noindent\underline{Proof of the approximation lemma \ref{lemme_approximation}}: 
Up to a negligible set, we can write
$$\r=\bigsqcup_{j\in\Z}I_j,$$
where $I_j=]a_j,a_{j+1}[$ is a bounded interval, $n^* (x)=n_j^*$ for $x\in
I_j$, and 
$n_j^*\in \{1,2\}$ (the assumption $n^*$ constant at $\pm\infty$ implies that
the sequence $(n_j^*)_{j\in\Z}$ is stationnary).\\
\noindent For all $k\in\N^*$, we can divide (up to a negligible set) each interval $I_j$
like this:
$$I_j=\bigsqcup_{i=0}^{k-1}]a_{j,i}^{(k)},a_{j,i+1}^{(k)}[,\qquad
a_{j,i}^{(k)}=a_j+\frac{i}{k}(a_{j+1}-a_j),\qquad i=0,\dots, k.$$
For $j\in\Z$, $k\in\N^*$, and $0\leq i\leq k-1,$ we set
$$
m_{j,i}^{(k)}=\frac{1}{n_j^*}\int_{a_{j,i}^{(k)}}^{a_{j,i+1}^{(k)}}n^0(x)dx.$$
Since $0\leq n^0\leq n^*$, we have $0\leq m_{j,i}^{(k)}\leq
\frac{meas(I_j)}{k}$. thus
\[]a_{j,i}^{(k)},a_{j,i}^{(k)}+m_{j,i}^{(k)}[\,\,\subset\,\,
]a_{j,i}^{(k)},a_{j,i+1}^{(k)}[.\]
We set 
$$n_k^0(x)=\sum_{j=-k}^{k}
\sum_{i=0}^{k-1}n_j^*\1_{]a_{j,i}^{(k)},a_{j,i}^{(k)}+m_{j,i}^{(k)}[}(x)
.$$
Obviously  $n_k^0$ satisfies $\eqref{approx_1}$.\\
\noindent
Moreover, we can notice that 
 $$\mbox{$n_k^0\equiv 0$ $a.e.$ on $]a_{j,i}^{(k)},a_{j,i+1}^{(k)}[$ $\ssi$
 $n^0\equiv 0$ $a.e.$ on $]a_{j,i}^{(k)},a_{j,i+1}^{(k)}[$},$$
and
 $$\mbox{$n_k^0\equiv n_j^*$ $a.e.$ on $]a_{j,i}^{(k)},a_{j,i+1}^{(k)}[$ $\ssi$
 $n^0\equiv n_j^*$ $a.e.$ on $]a_{j,i}^{(k)},a_{j,i+1}^{(k)}[$}.$$
We also define
$$n_k^0(x)u_k^0(x)=\sum_{j=-k}^{k}\sum_{i=0}^{k-1}n_j^*u_{j,i}^{(k)}\1_{]a_{j,i
}^{(k)},a_{j,i}^{(k)}+m_{j,i}^{(k)}[}(x),$$
where $\ds u_{j,i}^{(k)}=\underset{]a_{j,i}^{(k)},a_{j,i+1}^{(k)}[}{essinf}
u^0$, which makes sense because $u^0\in BV(\r)$.
We have 
 $$a.e.\,\, x\in ]a_{j,i}^{(k)},a_{j,i+1}^{(k)}[,\qquad \mbox{$n_k^0(x)\neq 0$
$\quad\Longrightarrow\quad$
 $u^0_k(x)=u_{j,i}^{(k)}$}.$$
We extend $u_k^0$ linearly in the vacuum (areas where $n_k^0= 0$) and at
infinity, as in Remark
\ref{prolongement_u_p}.\\
\noindent Thus, areas where $n_k^0=0$ have no influence on the total variation and we have 
$$\begin{array}[]{llll}
TV(u_k^0) & = & & |u_{-k,0}^{(k)}-u_{-k,1}^{(k)}|+\dots
 +|u_{-k,k-2}^{(k)}-u_{-k,k-1}^{(k)}|\\\\
& & + &|u_{-k,k-1}^{(k)}-u_{-k+1,0}^{(k)}|\\\\
& & +
&|u_{-k+1,0}^{(k)}-u_{-k+1,1}^{(k)}|+\dots+|u_{k-1,k-1}^{(k)}-u_{k,0}^{(k)}
|\\\\&  & + & |u_{k,0}^{(k)}-u_{k,1}^{(k)}|+\dots+
|u_{k,k-2}^{(k)}-u_{k,k-1}^{(k)}|\\\\
& \leq & & TV_{[a_{-k},a_{k+1}]}(u^0),
  \end{array}$$
which shows that $u_k^0$ satisfies $\eqref{approx_3}$. We also have
$\eqref{approx_2}$.\\
\noindent 
For any test function $\varphi\in\mathcal{D}(\r)$, we have
$$\begin{array}[]{lllll}
\ds\int_{\r}n_k^0(x)\varphi(x)dx & = &
\ds\sum_{|j|\leq k}\sum_{i=0}^{k-1}n_j^*\int_{a_{j,i}^{(k)}}^{a_{j,i}^{(k)}+
m_{j,i}^{(k)}}
\varphi(x)dx\\\\
& = & \ds\sum_{|j|\leq
k}\sum_{i=0}^{k-1}n_j^*\left(m_{j,i}^{(k)}\,\varphi(a_{j,i}^{(k)})
+\frac{{m_{j,i}^{(k)}}^2}{2}\,\varphi'({\xi_{j,i}^{(k)}})\right)
  \end{array}$$
with $a_{j,i}^{(k)}<{\xi_{j,i}^{(k)}}<a_{j,i}^{(k)}+m_{j,i}^{(k)}$ (if
${m_{j,i}^{(k)}}\neq 0$).
Thus, we can rewrite
$$\begin{array}[]{llll}
\ds\int_{\r}n_k^0(x)\varphi(x)dx & = & \ds\sum_{|j|\leq
k}\sum_{i=0}^{k-1}\left(\int_{a_{j,i}^{(k)}}^{a_{j,i+1}^{(k)}}n^0(x)
\varphi(a_{j,i}^{(k)}
)dx+\frac{n_j^*{m_{j,i}^{(k)}}^2}{2}\,\varphi'({\xi_{j,i}^{(k)}})\right).
  \end{array}$$
Let $j_0\in\N^*$ such that $\ds
supp(\varphi)\subset\bigsqcup_{|j|\leq j_0}I_{j_0}$
(it is possible because $\ds\inf_{j\in\Z}
\left(meas(I_j)\right)>0$).
Then
we have, for all $k\geq j_0$,
$$\begin{array}[]{llll}
\ds\int_{\r}n_k^0(x)\varphi(x)dx & = & \ds\sum_{|j|\leq
j_0}\sum_{i=0}^{k-1}\left(\int_{a_{j,i}^{(k)}}^{a_{j,i+1}^{(k)}}n^0(x)
\varphi(a_{j,i}^{(k)}
)dx+\frac{n_j^*{m_{j,i}^{(k)}}^2}{2}\,\varphi'({\xi_{j,i}^{(k)}})\right).
  \end{array}$$
We also have
$$\begin{array}[]{llll}
\ds\int_{\r}n^0(x)\varphi(x)dx & = & \ds\sum_{|j|\leq
j_0}\sum_{i=0}^{k-1}\int_{a_{j,i}^{(k)}}^{a_{j,i+1}^{(k)}}
n^0(x)\varphi(x)dx.
  \end{array}$$
Thus
$$\begin{array}[]{llll}
& & \ds\left|\int_{\r}n^0(x)\varphi(x)dx-\int_{\r}n_k^0(x)\varphi(x)dx\right|\\\\ 
& \leq
& \ds\sum_{|j|\leq
j_0}\sum_{i=0}^{k-1}\int_{a_{j,i}^{(k)}}^{a_{j,i+1}^{(k)}}n^0(x)
|\varphi(x)-\varphi(a_{j,i}^{(k)})|dx +   \,\,\,\,\ds \|\varphi'\|_{\infty}\ds\sum_{|j|\leq
j_0}\sum_{i=0}^{k-1}\frac{n_j^*{m_{j,i}^{(k)}}^2}{2}\\\\
&  \leq
&  \ds \|\varphi'\|_{\infty}\ds\ds\sum_{|j|\leq
j_0}\sum_{i=0}^{k-1}n_j^*\int_{a_{j,i}^{(k)}}^{a_{j,i+1}^{(k)}}(x-a_{j,i}^{(k)})
dx +   \,\,\,\,\ds \|\varphi'\|_{\infty}\ds\sum_{|j|\leq
j_0}\sum_{i=0}^{k-1}{m_{j,i}^{(k)}}^2\\\\
& \leq  & \ds 2\|\varphi'\|_{\infty} \ds\sum_{|j|\leq
j_0}\sum_{i=0}^{k-1}
\left(\frac{meas(I_j)}{k}\right)^2\\\\
&\ds \leq & C(\varphi,j_0)\times\ds\frac{1}{k}.
  \end{array}$$
Moreover, we have similarly
$$\begin{array}[]{llll}
& & \hspace{-15mm}\ds\int_{\r}n_k^0(x)u_k^0(x)\varphi(x)dx\\\\ 
& = & \ds\sum_{|j|\leq
j_0}\sum_{i=0}^{k-1}\left(\int_{a_{j,i}^{(k)}}^{a_{j,i+1}^{(k)}}n^0(x)
{u_{j,i}^{(k)}}
\varphi(a_{j,i}^{(k)}
)dx+\frac{n_j^*{m_{j,i}^{(k)}}^2}{2}{u_{j,i}^{(k)}}\,\varphi'({\xi_{j,i}^{(k)}}
)\right)
  \end{array}$$
and 
$$\begin{array}[]{llll}
\ds\int_{\r}n^0(x)u^0(x)\varphi(x)dx & = & \ds\sum_{|j|\leq
j_0}\sum_{i=0}^{k-1}\int_{a_{j,i}^{(k)}}^{a_{j,i+1}^{(k)}}
n^0(x)u^0(x)\varphi(x)dx
  \end{array}.$$
Thus
$$\begin{array}[]{llll}
& & \hspace{-15mm}\ds\left|\int_{\r}n^0(x)u^0(x)\varphi(x)dx-\int_{\r}
n_k^0(x)u_k^0(x)\varphi(x)dx\right| \\\\
&  \leq
& \ds\sum_{|j|\leq
j_0}\sum_{i=0}^{k-1}\int_{a_{j,i}^{(k)}}^{a_{j,i+1}^{(k)}}n^0(x){u_{j,i}^{(k)}}
|\varphi(x)-\varphi(a_{j,i}^{(k)})|dx\\\\
&  & +\quad\ds \|\varphi'\|_{\infty}\ds\sum_{|j|\leq
j_0}\sum_{i=0}^{k-1}{u_{j,i}^{(k)}}\frac{n_j^*{m_{j,i}^{(k)}}^2}{2}\\\\
& & +\quad \ds\sum_{|j|\leq
j_0}\sum_{i=0}^{k-1}\int_{a_{j,i}^{(k)}}^{a_{j,i+1}^{(k)}}n^0(x)|u^0(x)-
u_{j,i}^{(k)}||\varphi(x)|dx\\\\
& \leq & \ds C(\varphi,j_0)\times \|u_0\|_{\infty}\times \frac{1}{k}
\,\,\,\,+\,\,\,\,\ds2\|\varphi\|_{\infty}\sum_{|j|\leq
j_0}\sum_{i=0}^{k-1}\int_{a_{j,i}^{(k)}}^{a_{j,i+1}^{(k)}}|u^0(x)-
u_{j,i}^{(k)}|dx.
  \end{array}$$
Therefore we just need to show that the last term vanishes when
$k\rightarrow\infty$. This is raised because
$$\begin{array}[]{llll}
& &\hspace{-20mm} \ds\sum_{|j|\leq
j_0}\sum_{i=0}^{k-1}\left(\int_{a_{j,i}^{(k)}}^{a_{j,i+1}^{(k)}} |u^0(x)-
u_{j,i}^{(k)}|dx\right)\\\\
& \leq & \ds\sum_{|j|\leq
j_0}\sum_{i=0}^{k-1}\int_{a_{j,i}^{(k)}}^{a_{j,i+1}^{(k)}}
\left|\sup_{]a_{j,i}^{(k)},a_{j,i+1}^{(k)}[}
u^0-\inf_{]a_{j,i}^{(k)},a_{j,i+1}^{(k)}[} u^0\right|dx\\\\
& \leq & \ds\sum_{|j|\leq
j_0}\frac{meas(I_j)}{k}\left(\sum_{i=0}^{k-1}TV_{]a_{j,i}^{(k)},
a_{j,i+1}^{(k)}[}\left(u^0\right)\right)\\\\
& \leq & \ds\sum_{|j|\leq
j_0}\frac{meas(I_j)}{k}\,\,TV_{I_j}\left(u^0\right)\\\\
& \leq & \ds TV(u^0)\times C(j_0)\times \frac{1}{k}.
\end{array}$$
We established that $<n_k^0,\varphi>\,\,\rightarrow\,\, <n^0,\varphi>$ and
$<n_k^0u_k^0,\varphi>\,\,\rightarrow\,\, <n^0u^0,\varphi>$.
Finally, we define $p_k^0$  the same way as $u_k^0$:
$$n_k^0(x)p_k^0(x)=\sum_{j=-k}^{k}\sum_{i=0}^{k-1}n_j^*p_{j,i}^{(k)}\1_{]a_{j,i
}^{(k)},a_{j,i}^{(k)}+m_{j,i}^{(k)}[}(x),$$
where $\ds p_{j,i}^{(k)}=\underset{]a_{j,i}^{(k)},a_{j,i+1}^{(k)}[}{essinf}
p^0$. But in the vacuum (areas where $n_k^0=0$) we set $p_k^0=0.$
Thus, we have $p_k^0\equiv p_{j,i}^{(k)}$ on each interval
$]a_{j,i}^{(k)},a_{j,i+1}^{(k)}[$. In fact, there are two cases:
\begin{itemize}
 \item If $n^0\equiv n_j^*$ a.e. on $]a_{j,i}^{(k)},a_{j,i+1}^{(k)}[$, then
$n_k^0\equiv n_j^*$ and $p_k^0 \equiv p_{j,i}^{(k)}.$
\item Else, it exists a non negligible subset $\omega\subset
]a_{j,i}^{(k)},a_{j,i+1}^{(k)}[$ where $n^0<n_j^*$, and $p^0\equiv 0$ a.e. on
$\omega$, which implies $p_{j,i}^{(k)}=0$, and 
$p_k^0\equiv 0=p_{j,i}^{(k)}$ a.e. on $]a_{j,i}^{(k)},a_{j,i+1}^{(k)}[$.
\end{itemize}
We easily deduce that $p_k^0$ satisfies properties $\eqref{approx_2}$ and
$\eqref{approx_3}$.

For the
convergence $<n_k^0p_k^0,\varphi>\,\,\rightarrow\,\, <n^0p^0,\varphi>$, the
proof is
exactly the same as $n_k^0u_k^0$. Finally, the last point is obvious because 
$n_k^0(x)\in\{0,n_j^*\}$ for all $x\in\r$, thus we have 
$$(n^*(x)-n_k^0)p_k^0=0,\qquad\forall k\geq 1.$$
\cqfd
 
\noindent\underline{Proof of the lemma \ref{lemme de compacite}}: Let $(\varphi_m)_{m\geq 1}$ be a countable set dense in
$\mathcal{D}(\r_x)$ for the $L^1$-norm, which exists because of the
separability of $L^1(\r_x)$. We denote
$$g_{k,m}(t):=\int_{\r}n_k(t,x)\varphi_m(x)dx.$$
The sequence $(g_{k,1})_{k\geq 1}$ is bounded and
equicontinuous in $C([0,T],\r)$, thus, the Ascoli Theorem entails that it
exists an extraction $\sigma_1(k)$ such that 
$$g_{\sigma_1(k),1}\underset{k}{\longrightarrow} l_1\quad\mbox{in}\quad
C([0,T],\r).$$
The same applies to $(g_{\sigma_1(k),2})_{k\geq 1}$, thus it exists an
extraction $\sigma_2$ such that 
$$g_{\sigma_1(\sigma_2(k)),2}\underset{k}{\longrightarrow}
l_2\quad\mbox{in}\quad
C([0,T],\r).$$
A simple recursion shows that we can build a sequence of extractions $\sigma_m$
such that 
$$g_{\sigma_1(\sigma_2(\dots\sigma_m(k))\dots)),m}\underset{k}{\longrightarrow}
l_m\quad\mbox{in}\quad
C([0,T],\r).$$
Therefore, setting $\sigma(k):=\sigma_1\circ\dots\circ\sigma_k(k),$ we have (by
diagonal extraction)
\begin{equation}\label{cv unif}
\forall m\geq 1,\qquad g_{\sigma(k),m}\underset{k}{\longrightarrow}
l_m\quad\mbox{in}\quad
C([0,T],\r). 
\end{equation}
Now, we can identify the limit $l_m$ because since
$(n_{\sigma(k)})_k$ is bounded in $L^{\infty}(]0,T[\times \r)$, there exists a
subsequence (still denoted by the same way) such that $n_{\sigma(k)}\tw n$ in
$L^{\infty}_{w*}(]0,T[\times \r)$. Thus, we have, for all $m\geq 1$, and for
all $\psi\in\mathcal{D}(]0,T[_t)$,
$$\int_{0}^{T}\int_{\r}n_{\sigma(k)}(t,x)\psi(t)\varphi_m(x)dxdt\underset{
k}{\rightarrow}\int_{0}^{T}\int_{\r}n(t,x)\psi(t)\varphi_m(x)dxdt,$$
which rewrites
$$\int_{0}^T g_{\sigma(k),m}(t)\psi(t)dt\underset{k}{\rightarrow}\int_0^T
\left(\int_{\r}n(t,x)\varphi_m(x)dx\right)\psi(t)dt.$$
Moreover, \eqref{cv unif} easily implies that 
$$\int_0^T g_{\sigma(k),m}(t)\psi(t)dt\underset{k}{\rightarrow}\int_0^T
l_m(t)\psi(t)dt,$$
thus $l_m(t)=\int_{\r}n(t,x)\varphi_m(x)dx$, a.e. $t\in[0,T]$, from which  we can
deduce
$$\forall m\geq 1,\qquad \sup_{t\in[0,T]}\left|\int_{\r}(n_{\sigma(k)}(t,x)-n(t,
x))\varphi_m(x)dx\right|\underset{k}{\longrightarrow} 0.$$
Finally, this convergence stays available for all
$\varphi\in\mathcal{D}(\r_x)$, because of the inequality
$$\begin{array}{lll}
\ds\sup_{t\in[0,T]}\left|\int_{\r}(n_{\sigma(k)}-n)(t,x)\varphi(x)dx\right|\\\\ 
\leq
\ds\sup_{t\in[0,T]}\left|\int_{\r}(n_{\sigma(k)}
-n)(t,x)\varphi_m(x)dx\right|+C\|\varphi-\varphi_m\|_{L^1(\r)},
\end{array}$$
where \[C:=\sup_{k\geq 1}(\|n_k\|_{L^{\infty}(]0,T[\times
\r)})+\|n\|_{L^{\infty}(]0,T[\times \r)}<+\infty.\]
We conclude that it is also true for $\Gamma\in L^1(\r_x)$ by density, using
the same inequality. 
\cqfd

\end{document}